\documentclass[12pt,reqno]{amsart}

\usepackage{amssymb,amsmath,amsbsy,eucal,amsfonts,mathrsfs,latexsym}
\usepackage{graphicx,verbatim,psfrag,epsfig,enumerate}
\usepackage{stmaryrd}
\usepackage{lscape}
\usepackage{url}
\usepackage{multirow}
\usepackage{color}
\usepackage{amsmath,amssymb,amsfonts}
\usepackage{graphicx}
\usepackage{subcaption}
\usepackage{booktabs}
\usepackage{array}
\usepackage{longtable}
\usepackage{multirow}
\usepackage{siunitx}
\usepackage{xcolor}
\definecolor{reviewblue}{RGB}{0,70,160}

\usepackage{float}
\usepackage{placeins}
\usepackage{enumitem}
\usepackage{longtable}
\usepackage{booktabs}
\usepackage{subcaption}
\graphicspath{{graph/}}

\usepackage{booktabs}
\usepackage{tabularx}
\usepackage{adjustbox}

\newtheorem{theorem}{Theorem}[section]

\theoremstyle{definition}

\numberwithin{equation}{section}

\newcommand{\vertiii}[1]{{\left\vert\kern-0.25ex\left\vert\kern-0.25ex\left\vert #1 
    \right\vert\kern-0.25ex\right\vert\kern-0.25ex\right\vert}}

\definecolor{red}{rgb}{1,0,0}

\definecolor{blue}{rgb}{0,0,1}

\usepackage{mathtools,xcolor}
\usepackage[colorlinks=true,linkcolor=blue]{hyperref}
\usepackage{cleveref}
\usepackage{algorithm}
\usepackage{algorithmic}

\begin{document}

\title{A Neural-network-based multiscale Hybridizable Discontinuous Galerkin method for solving PDEs in porous media}

\author{Tony Haines}
\address{Department of Mathematics, ECPI University}
\email{thaines@ecpi.edu}

\author{Ke Shi}
\address{Department of Mathematics $\&$ Statistics, Old Dominion University, Norfolk, VA 23529, USA}
\email{kshi@odu.edu}
\thanks{
 As a convention the names of the authors are alphabetically ordered. 
Both authors contributed equally in this article. 
}


\begin{abstract}
We develop a neural-network-accelerated multiscale hybridizable
discontinuous Galerkin method for elliptic problems with heterogeneous
coefficients.  The method preserves the standard MsHDG local-to-global structure: fine-scale HDG problems on coarse blocks define discrete Dirichlet-to-Neumann operators, which are assembled through the standard MsHDG global skeleton equations. To reduce the cost of constructing these local
operators, we train a neural network on coefficient fields defined on a
reference block and use the predicted operators in place of repeated
fine-scale local solves.

The numerical experiments assess both the accuracy and online efficiency of the resulting NN-MsHDG method. For two-dimensional binary permeability fields with moderate contrast, the neural method reproduces the standard MsHDG solution with relative modeling errors of a few percent while reducing the total online computational cost by factors of approximately 5 to 16, depending on the coarse trace dimension. In the high-contrast regime, the chosen polynomial coarse trace spaces already yield substantial discretization errors in the standard MsHDG method, indicating the need for more effective coarse spaces, such as coefficient-adapted spectral trace spaces. In addition, the learned local operators introduce modeling errors that become severe as the trace space is enriched. These results demonstrate the potential of neural surrogates for accelerating multiscale HDG computations while also highlighting the need for improved operator representations and a better understanding of error amplification in high-contrast problems.
\end{abstract}

\subjclass[2020]{65N30, 65N55, 68T07, 76S05}

\keywords{Multiscale HDG, hybridizable discontinuous Galerkin methods,
neural networks, Dirichlet-to-Neumann maps, heterogeneous media.}

\maketitle

\section{Introduction}

Many physical models are governed by partial differential equations whose
coefficients contain strong spatial heterogeneity.  Examples include flow in
porous media, heat conduction in composite materials, diffusion and transport
in heterogeneous media, and wave propagation through materials with complicated
microstructure.  In such problems, the coefficient field may vary over a wide
range of spatial scales and may also exhibit large contrast between different
regions of the domain.  A fine-scale numerical discretization that resolves all
relevant heterogeneities can therefore lead to very large algebraic systems.
This difficulty becomes even more severe in parameter-dependent,
time-dependent, or nonlinear problems, where related local or global solves may
have to be repeated many times.

Hybridizable discontinuous Galerkin (HDG) methods provide a natural
framework for reducing the number of globally coupled degrees of
freedom; see, for example,
\cite{CockburnGopalakrishnanLazarov2009,
CockburnGopalakrishnanSayas2010,
CockburnQiuShi2012}.
In an HDG method, the elementwise or subdomainwise unknowns are
eliminated locally, and the global system is written only in terms of
trace unknowns on the mesh skeleton.  Once the trace is determined, the
remaining unknowns can be recovered independently by local solvers.
This hybridized structure may be interpreted as a discrete
Dirichlet-to-Neumann or input-output formulation: for a prescribed trace
on the boundary of a local region, one solves a local problem and
computes the corresponding numerical flux or output response.  The
global problem is then obtained by enforcing conservation of these local
fluxes across interfaces.  The HDG framework has also been developed for systems and nonlinear
problems, including the Stokes and Maxwell equations, nonlinear
convection--diffusion equations, the incompressible Navier--Stokes
equations, and nonlinear elasticity; see, for example,
\cite{NguyenPeraireCockburn2010,
NguyenPeraireCockburn2011Maxwell,
NguyenPeraireCockburn2009Nonlinear,
NguyenPeraireCockburn2011NavierStokes,
KabariaLewCockburn2015}.

For problems in heterogeneous media, multiscale HDG methods further
exploit this local-to-global structure
\cite{EfendievLazarovShi2015,EfendievLazarovMoonShi2014,
YangShiFu2019,MaFu2021}.  The original MsHDG framework uses polynomial
or homogenization-based coarse trace spaces
\cite{EfendievLazarovShi2015}, while coefficient-adapted spectral trace
spaces were subsequently developed for high-contrast media
\cite{EfendievLazarovMoonShi2014}.  Related extensions have been
introduced for Darcy and two-phase flow
\cite{YangShiFu2019} and for heterogeneous linear elasticity
\cite{MaFu2021}.  In these methods, the domain is decomposed into coarse
blocks, each equipped with a fine mesh capable of resolving the local
coefficient variation.  The globally coupled unknown is a coarse trace
variable defined on the coarse skeleton.  For each prescribed coarse
trace, a fine-scale HDG problem is solved inside each coarse block, and
the resulting normal numerical flux defines a coarse-block
Dirichlet-to-Neumann map.  Thus, the MsHDG method reduces the global
problem to a system posed on the coarse skeleton while retaining
fine-scale information through the local solvers.

The computational bottleneck of this procedure lies in the construction of the
local coarse-block solution operators.  Even though the final global system has
far fewer degrees of freedom than the fine-scale system, assembling the coarse
operator requires solving many fine-scale local problems.  In particular, for
each coarse block one must apply the inverse of a fine-scale local HDG matrix to
all coarse trace basis inputs.  This cost must be paid separately for each local
coefficient configuration.  In stationary linear problems with fixed
coefficients this local construction is performed once.  However, in
time-dependent or nonlinear problems the local operators may change at each time
step or nonlinear iteration.  For example, the local operator may depend on an
evolving coefficient, on the solution from the previous time step, or on the
current linearization state.  In such settings, repeated local-solver
construction can dominate the total computational cost.

There has been rapidly growing interest in using neural networks to accelerate
or approximate the numerical solution of partial differential equations.  One
line of work uses neural networks directly as trial functions for the solution,
leading to methods such as physics-informed neural networks
\cite{RaissiPerdikarisKarniadakis2019}, the Deep Ritz method
\cite{Yu2018}, and deep learning methods for high-dimensional parabolic
equations \cite{EHanJentzen2017,HanJentzenE2018}.  These methods are attractive
because of the flexibility of neural-network approximation spaces, but their
performance may depend strongly on the optimization landscape, the choice of
loss function, and the enforcement of boundary conditions or conservation laws.
A second line of work aims to learn solution operators rather than individual
solutions.  This includes DeepONet \cite{LuJinPangZhangKarniadakis2021},
neural operators \cite{KovachkiLiLiuAzizzadenesheliBhattacharyaStuartAnandkumar2023},
Fourier neural operators
\cite{LiKovachkiAzizzadenesheliLiuBhattacharyaStuartAnandkumar2021}, and
geometry-aware operator learning methods
\cite{LiKovachkiChoyEtAl2023,SerranoLeBoudecKassaiKoupaiWangYinVittautGallinari2023}.
Once trained, such methods can provide fast online evaluation of
parameter-to-solution maps, although their generalization may be limited by the
training distribution, geometry, and boundary conditions.

A third direction, closer to the present work, seeks to combine neural networks
with established numerical discretizations rather than replace them entirely.
Examples include neural-network solvers for parametric PDEs
\cite{KhooLuYing2021}, machine-learning-accelerated computational fluid
dynamics \cite{KochkovSmithAlievaWangBrennerHoyer2021}, neural solvers for
radiative transfer \cite{TanoRagusa2021,LuWangXu2022}, and learned local or
reduced solvers for multiscale problems \cite{ChenDingLiWright2024}.  The
advantage of this viewpoint is that one can retain the structure, geometry handling, and coupling mechanisms of classical numerical
methods, while using neural networks only to accelerate selected expensive
components.

The work of Du and Stechmann \cite{DuStechmann2024} follows this hybrid
philosophy.  They proposed an element-learning framework in which a neural
network is trained on a reference element to predict local input-output
operators.  For transport-type problems these include the inflow-to-outflow
operator and the inflow-to-solution operator; for elliptic-type problems, the
corresponding input-output operator can be viewed as a Dirichlet-to-Neumann
map.  Their method is closely related to standard HDG methods, since the neural
network replaces the construction of local HDG solvers on each fine-scale
element.  Their numerical experiments for radiative transfer demonstrate that
such a strategy can provide substantial speedup when the local element problems
are sufficiently expensive.

The present work is motivated by the same general principle--replacing
expensive hybridized local solution operators by neural-network surrogates--but
applies it at a different level of the discretization.  In the element-learning
approach of \cite{DuStechmann2024}, the neural network approximates the local
solution operator on each fine-scale element of a standard HDG method.  This is
especially advantageous when each element contains many local degrees of
freedom, for example in high-order spectral elements or in kinetic and
radiative-transfer models with additional angular variables.  In contrast, many
multiscale elliptic problems in heterogeneous media are effectively discretized
by low-order finite elements on a fine mesh.  In this regime, the local solve
on a single fine element may be too small for a neural-network replacement to
offer a significant advantage.  The expensive object is instead the local
solution operator on an entire coarse block, which contains many fine elements
and encodes the effect of the heterogeneous medium on the coarse skeleton.

Another work closely related to the present paper is
\cite{BoutilierBrennerMiguez2025}, where neural networks are used to learn
local nonlinear Dirichlet-to-Neumann maps for elliptic problems with rough
coefficients.  Their method is formulated in a multiscale finite element
substructuring framework.  The local problems are discretized by conforming
finite elements, and the learned DtN maps are used in a nonlinear substructured
problem.  Our work shares the same local-to-global philosophy: expensive local
DtN evaluations are replaced by neural-network surrogates.  However, the local
discretization and the learned operator are different.  We work within the
multiscale HDG framework, where the local fine-scale unknowns inside each
coarse block are eliminated by HDG local solvers and the global unknown is a
hybrid trace on the coarse skeleton.  For the linear Darcy model considered
here, the local trace-to-flux map is represented by a discrete coarse-block DtN
matrix \(S_i(\alpha)\).  Thus the neural network is trained to approximate the
coefficient-to-operator map
\[
    \alpha|_{\Omega_i} \longmapsto S_i(\alpha),
\]
rather than a conforming finite element nonlinear trace-to-flux map.  The
resulting learned operators are assembled directly into the standard MsHDG
coarse skeleton system.

We therefore propose a neural-network-based multiscale HDG method
(NN-MsHDG).  For each coarse block, the fine-scale HDG local solver
produces a discrete DtN matrix \(S_i(\alpha)\), describing the
trace-induced flux response, together with a vector
\(g_i(\alpha,f)\), describing the source-induced contribution.  The
input to the network is a representation of the local coefficient
field, and the outputs are the independent entries of
\(S_i(\alpha)\) and the components of \(g_i(\alpha,f)\).  Once trained
on representative local coefficient configurations, the network can
rapidly generate the quantities required for the global skeleton
assembly.  The predicted local operators are assembled exactly as in
the standard MsHDG method, so the global hybridized coupling and local
reconstruction procedures remain unchanged.

The distinction between element-level learning, conforming finite element
substructuring, and coarse-block multiscale HDG learning is central to the
present approach.  In standard element-level learning, the neural network
replaces many small HDG element solvers, and the computational advantage is
most visible when the local element dimension is large.  In conforming finite
element substructuring approaches, neural networks may approximate local DtN
maps associated with nodal trace variables and conforming local solves.  In the
present MsHDG setting, the neural network replaces a coarse-block DtN
construction, where each local map is obtained from a fine-scale HDG
discretization over an entire subdomain.  Consequently, the approach remains
meaningful even when the fine-scale HDG spaces are low order, because the local
coarse-block operator is still expensive to construct.  This makes the proposed
framework particularly suitable for multiscale heterogeneous problems where
fine-scale accuracy is needed locally but only coarse trace unknowns are
globally coupled.

Although the numerical examples in this paper focus on stationary Darcy flow
for clarity, the framework is not restricted to this model problem.  Darcy flow
serves as a clean prototype because its HDG local solver naturally defines a
Dirichlet-to-Neumann map and because heterogeneous permeability fields provide a
standard multiscale setting.  The potential advantage of the method is expected
to be even greater for time-dependent or nonlinear problems in which the local
operators must be updated repeatedly.  Examples include nonlinear diffusion,
flow with state-dependent permeability, multiphase porous media flow after
linearization, reactive transport, and heat conduction with temperature-dependent
conductivity.  In such problems, a neural network can be trained to approximate
the parameter-to-operator map associated with the local HDG solver, thereby
turning repeated local matrix inversions into inexpensive online evaluations.

The remainder of the paper is organized as follows.
Section~\ref{sec:general-msHDG} reviews the multiscale HDG framework
and its local-to-global structure.  Section~\ref{sec:darcy-flow}
specializes the formulation to the Darcy flow problem and introduces
the corresponding fine- and coarse-scale HDG spaces.
Section~\ref{sec:mutiscaleHDG} derives the coarse-block local system in
matrix form and identifies the discrete DtN matrix and source vector
required for the global skeleton assembly.  Section \ref{sec:NN solvers} describes the
reference-block data generation, neural-network architecture, and
training procedure used to approximate these local operators.
Section \ref{sec:numerical-experiments} presents the numerical experiments, including the
one-dimensional fixed-sample coefficient-resolution study, the two-dimensional lower-contrast
accuracy and timing studies, and the high-contrast tests of both local
operator action and assembled solution accuracy.  Section \ref{sec:conclusion} summarizes
the main findings and discusses possible extensions and remaining
challenges, particularly for high-contrast media.

\section{Review of the MsHDG methods for PDEs in heterogeneous media}\label{sec:general-msHDG}

We first describe the MsHDG framework at the level of a general
steady-state boundary-value problem.  Let
\(\Omega\subset\mathbb R^d\), \(d=2,3\), be a bounded polygonal or
polyhedral domain, and consider
\[
    \mathcal L_\kappa u=f
    \quad\text{in }\Omega,
    \qquad
    u=g
    \quad\text{on }\partial\Omega,
\]
where \(\kappa\) denotes the possibly heterogeneous physical
coefficients.  This notation directly covers linear and nonlinear elliptic and
diffusion problems.  The same local-to-global construction also extends
to systems, such as the Navier-Stokes equations, whenever the corresponding
HDG local solver defines an appropriate trace-to-flux or
trace-to-traction map.  At this stage, we do not specify the weak
formulation or finite element spaces. Instead, we emphasize the domain decomposition idea that underlies the MsHDG methodology. The subscript \(\kappa\) collectively denotes the possibly
heterogeneous physical parameters entering the differential operator.

\begin{figure}[ht] 
 \includegraphics[ width=.4\textwidth, clip]{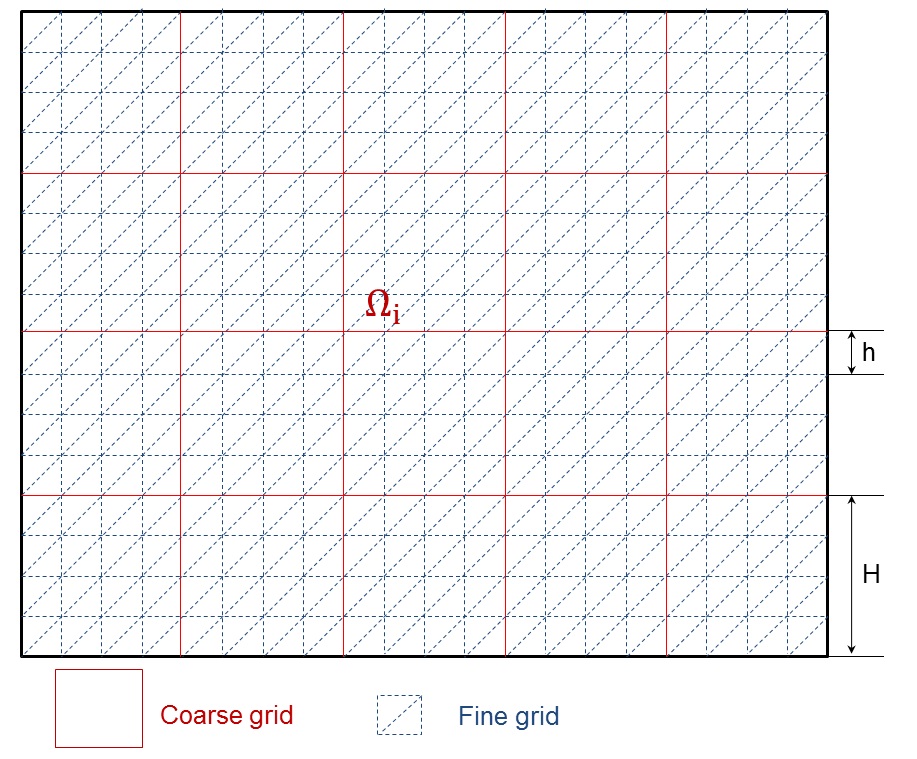}
\centering
\vskip5mm
\caption{Partition of $\Omega$ with coarse blocks $\Omega_i$. Each $\Omega_i$ is partitioned with fine scale triangulations.}
\label{fig:partition}
\end{figure}

Let
$
    \mathcal P=\{\Omega_i\}_{i=1}^{N_c}
$
be a nonoverlapping coarse partition of \(\Omega\), so that
$
    \overline{\Omega}
      =\bigcup_{i=1}^{N_c}\overline{\Omega_i}.
$
The union of all coarse interfaces is denoted by \(\mathcal E\) and is
called the coarse skeleton.  Each interface \(F\in\mathcal E\) is
partitioned into trace faces \(F_H\), and the collection of all such
faces is denoted by \(\mathcal E_H\)   as illustrated in Figure~\ref{fig:partition}.

 The main global unknown in the MsHDG framework is a trace variable
\[
    \mu_H \in M_H,
\]
where $M_H$ is an admissible finite dimensional space with coarse scale H defined on $\mathcal{E}_H$. At the PDE level, this trace variable is treated as Dirichlet input data for the local subdomain problems.

More precisely, for each coarse subdomain $\Omega_i$ and for a prescribed trace value $\mu$ on $\partial \Omega_i$, we define the local problem
\begin{equation}\label{eq:abstract-local-problem}
    \mathcal{L}_{\kappa} u = f \quad \text{in } \; \Omega_i,
    \qquad
    u = \mu_H \quad \text{on } \partial \Omega_i.
\end{equation}
The solution of this local problem determines a Neumann response, or more generally a flux response, which we denote by $q_n$. The precise meaning of $q_n$ depends on the equation under consideration. For Darcy flow it is the normal Darcy flux; for the Stokes equations it is related to the normal stress; and for other second-order PDEs it represents the corresponding natural boundary quantity. Thus, the local solution operator defines a Dirichlet-to-Neumann (DtN) map
\[
    (\Omega_i, \kappa, f, \mu) \longmapsto q_n.
\]

Since the local problems are solved independently on the coarse subdomains, the flux $q_n$ is not automatically continuous across the interior coarse interfaces. The global coupling condition is imposed by seeking a trace $\mu_H \in M_H$ such that the corresponding flux is conservative across all interior interfaces, namely the {\em transmission equation}:
\begin{equation}\label{eq:abstract-flux-continuity}
    \llbracket q_n \rrbracket = 0
    \quad \text{on all interior faces of } \mathcal{E}_H,
\end{equation}
with the prescribed boundary condition incorporated on $\partial \Omega$. Here $\llbracket \cdot \rrbracket$ denotes the jump of the normal flux across a coarse interface. This flux-continuity condition provides the global coupling between the independently computed local solutions. At the continuous level, the trace of the exact solution satisfies the
transmission condition, and uniqueness of the original boundary-value
problem determines this trace uniquely.  In the multiscale discretization, the trace is sought in the
finite-dimensional coarse space \(M_H\), and the transmission condition
is imposed weakly by testing the coarse-interface fluxes against all
functions in \(M_H\).  The resulting trace is an approximation of the
exact solution trace whose quality depends on the approximation
properties of \(M_H\).

Therefore, the MsHDG framework can be viewed as a domain decomposition method in which the globally coupled unknown lives only on the coarse skeleton $\mathcal{E}_H$. The fine-scale variables inside each coarse subdomain are eliminated through local solvers. Once the trace $\mu$ is determined, the solution inside each coarse subdomain can be recovered independently.

In practical computations, the Dirichlet-to-Neumann map (DtN) is not evaluated at the continuous PDE level. Instead, each coarse subdomain is equipped with a fine mesh that resolves the local heterogeneity, and the local boundary value problems are discretized by an HDG method. For a linear problem, the local flux response can be decomposed into a
trace-induced and a source-induced contribution,
\[
    \operatorname{DtN}_{\Omega_i,\kappa}(\mu_H,f)
      =
    \operatorname{DtN}_{\Omega_i,\kappa}(\mu_H,0)
      +
    \operatorname{DtN}_{\Omega_i,\kappa}(0,f).
\]
At the discrete level, these two contributions give the local matrix
\(S_i(\alpha)\) and the source vector \(g_i(\alpha,f)\), respectively.
Their precise construction is derived in Section \ref{sec:mutiscaleHDG}.

The discrete DtN maps provide the local contributions required for the
assembly of the global coarse-scale system.  Although the resulting
global system contains only the coarse trace unknowns, its construction
is not necessarily inexpensive.  For each coarse block and each local
coefficient configuration, the fine-scale HDG problem must be solved
repeatedly for the prescribed coarse trace basis functions, together
with the source-induced local problem.  The resulting local flux
responses are then assembled into the global skeleton matrix and
right-hand side.  Consequently, the construction and assembly of the
local DtN operators can constitute the dominant part of the online
computational cost.  This behavior is confirmed by the timing results
in Section~\ref{sec:numerical-experiments}, where the local-operator
construction and global assembly account for most of the standard
MsHDG runtime and become increasingly expensive as the coarse trace
space is enriched.

This computational burden is expected to be even more significant for
more complicated PDEs.  For systems with several physical variables,
high-order or three-dimensional discretizations, the local HDG matrices
are substantially larger.  Moreover, for time-dependent or nonlinear
problems, the local operators may have to be reconstructed at every
time step or nonlinear iteration.  Replacing these repeated local
computations by an efficient surrogate therefore has the potential to
produce substantial computational savings.

In the next section, we specialize the framework to the Darcy flow
problem and recall the corresponding MsHDG discretization developed
in~\cite{EfendievLazarovShi2015}.  Section~\ref{sec:mutiscaleHDG}
then derives the local coarse-block operators in matrix form and makes
explicit the computational bottleneck targeted by the neural-network
approximation.

\section{Darcy flow and the multiscale partitions}\label{sec:darcy-flow}

To make the idea precise, we consider the following second-order Darcy flow problem in porous media, on which the original MsHDG method was initially developed:
\begin{equation}\label{eq:general}
    - \nabla \cdot( \kappa(x) \nabla u) = f(x), \quad x \in \Omega,
\end{equation}
with homogeneous Dirichlet boundary conditions $u = 0$ on $\partial \Omega$. Here $\kappa(x) \ge \kappa_0 >0$ is a possibly highly heterogeneous coefficient and $\Omega$ is a bounded polyhedral domain in $\mathbb{R}^d$, $d=2,3$. The methods presented here target applications of equation \eqref{eq:general} to flows in porous media, while the same framework also applies to diffusion, passive chemical transport, and heat transfer in heterogeneous media.

As a standard procedure of HDG methods, we begin by writing this problem in mixed form:
\begin{subequations}\label{eq:main}
\begin{alignat}{2}
\label{original equation-1}
\alpha \boldsymbol{q} + \nabla u &= 0 \qquad && \text{in $\Omega$,}\\
\label{original equation-2}
\nabla \cdot \boldsymbol{q}  &= f && \text{in $\Omega$,}\\
\label{boundary condition}
u &= 0 && \text{on $\partial \Omega$.}
\end{alignat}
\end{subequations}
Here $\alpha(x)=\kappa(x)^{-1}$ and $f \in L^2(\Omega)$.

Let $\mathcal{P} = \{\Omega_i\}_{i=1}^{N_c}$ be a disjoint polygonal partition of the domain $\Omega$ with $\Omega = \cup_{i=1}^{N_c} \Omega_i$ as shown in Figure \ref{fig:partition}. For simplicity, we illustrate the method using a uniform partition. In general, the partition may be quasi-uniform and may contain hanging nodes.

We call $F$ an interface of the partition $\mathcal{P}$ if $F$ is 
either shared by two neighboring subdomains, $\Omega_i$ and $\Omega_j$ in $\mathcal{P}$, 
$F = \overline{\Omega}_i \cap \overline{\Omega}_j$, or it is on the boundary $\partial \Omega$,
$F = \overline{\Omega}_i \cap \partial{\Omega}$. 
Let $\mathcal{E}$ denote the skeleton (or all interfaces) of the partition $\mathcal{P}$. 
For each interface $F \in \mathcal{E}$, let $\mathcal{E}_{H}(F)$ be a quasi-uniform partition of $F$ with maximum 
size (diameter) $H$. Set $\mathcal{E}_H = \cup_{F \in \mathcal{E}} \mathcal{E}_{H}(F)$.

For each coarse block \(\Omega_i\), let
\(\mathcal T_h(\Omega_i)\) be a conforming, shape-regular
triangulation of \(\Omega_i\), and set
$
    \mathcal T_h
      :=
    \bigcup_{i=1}^{N_c}\mathcal T_h(\Omega_i).
$
Let
$
    h:=\max_{K\in\mathcal T_h}\operatorname{diam}(K).
$
We denote by \(\mathcal E_h(\Omega_i)\) the set of fine faces in
\(\mathcal T_h(\Omega_i)\), and by
$
    \mathcal E_h^0(\Omega_i)
      :=
    \{F\in\mathcal E_h(\Omega_i):F\not\subset\partial\Omega_i\}
$
the interior fine faces of \(\Omega_i\) and $\mathcal{E}^0_h:= \cup_{i=1}^{N_c} \mathcal{E}^0_h(\Omega_i)$.


\[
\begin{aligned}
W_h
  &:=
  \{w\in L^2(\Omega):
      w|_K\in W(K),\ K\in\mathcal T_h\},\\
\boldsymbol V_h
  &:=
  \{\boldsymbol r\in [L^2(\Omega)]^d:
      \boldsymbol r|_K\in\boldsymbol V(K),\
      K\in\mathcal T_h\},\\
M_h^0
  &:=
  \{\mu\in L^2(\mathcal E_h^0):
      \mu|_F\in M_h(F),\
      F\in\mathcal E_h^0\},\\
M_H
  &:=
  \{\mu\in L^2(\mathcal E_H):
      \mu|_{F_H}\in M_H(F_H),\
      F_H\in\mathcal E_H,\
      \mu|_{\partial\Omega}=0\},\\
M_{h,H}
  &:=M_h^0\oplus M_H.
\end{aligned}
\]


The multiscale  HDG method is:
Find $(u_h, \boldsymbol{q}_h, \widehat{u}_{h,H}) \in W_h \times \boldsymbol{V}_h \times M_{h,H}$
such that
\begin{subequations}\label{eq:fineHDG}
\begin{alignat}{5}
\label{eq:fineHDG-1}
 &(\alpha \boldsymbol{q}_h , \boldsymbol{r})_{\mathcal{T}_h} - &&(u_h, \nabla \cdot \boldsymbol{r})_{\mathcal{T}_h} + &&\left<\widehat{u}_{h,H}, \boldsymbol{r} \cdot \boldsymbol{n}\right>_{\partial\mathcal{T}_h}    &&=0     && \forall \boldsymbol{r} \in \boldsymbol{V}_h,\\
 \label{eq:fineHDG-2}
-&(\boldsymbol{q}_h, \nabla w)_{\mathcal{T}_h}           &&                                       +&&\left<\widehat{\boldsymbol{q}}_{h,H} \cdot \boldsymbol{n}, w\right>_{\partial\mathcal{T}_h}    &&=(f, w)_{\mathcal{T}_h} && \forall w \in W_h,\\
\label{eq:fineHDG-3}
 &                                                                              &&  &&\left<\widehat{\boldsymbol{q}}_{h,H} \cdot \boldsymbol{n}, \mu\right>_{\partial\mathcal{T}_h}  &&=0     && \forall \mu \in M_{h,H},\\
 \label{eq:fineHDG-4}
 &                                                                              && && ~~\widehat{u}_{h,H}                                                      &&=0     && ~~\mbox{on}~~\partial \Omega.
\end{alignat}
\end{subequations}
Here,  $\widehat{u}_{h,H}$ and $ \widehat{\boldsymbol{q}}_{h,H} $ are the
{\it numerical trace} and  the {\it numerical flux}, respectively.
We write $(\eta, \zeta)_{\mathcal{T}_h}:= \sum_{K \in \mathcal{T}_h} (\eta, \zeta)_K,$ where $(\eta, \zeta)_D$ denotes the integral of $\eta\zeta$ over the domain $D \subset \mathbb{R}^d$
and
$\left<\eta, \zeta \right>_{\partial\mathcal{T}_h}
:=\sum_{K \in \mathcal{T}_h}\left<\eta, \zeta \right>_{\partial K}$,
where $\left<\eta, \zeta \right>_{\partial D}$ denotes
the integral of $\eta\zeta$ over the boundary of the domain $\partial D \subset \mathbb{R}^{d-1}$.
Consequently,
$\left<\eta, \zeta \right>_{\partial\Omega_i}:=\int_{\partial \Omega_i} \eta \zeta ds$.
The HDG method formulation is completed with the definition of the normal
component of the numerical flux :
\begin{equation}\label{eq:trace}
\widehat{\boldsymbol{q}}_{h,H} \cdot \boldsymbol{n} = \boldsymbol{q}_h \cdot \boldsymbol{n} + \tau (u_h - \widehat{u}_{h,H}),
\end{equation}
%
where $\tau$ is a non-negative stabilization parameter. More precisely, on each fine-scale edge $F \in \mathcal{E}_h$, $\tau$ is a fixed non-negative number. 
We remark that if there is only one coarse
block (i.e. $\Omega_i=\Omega$,  $i=1$) then this is the standard HDG
method \cite{CockburnQiuShi2012}. 

\subsection{Well-posedness and Convergence of the MsHDG}

The above scheme \eqref{eq:fineHDG}, \eqref{eq:trace} is complete with explicit choices of the stabilization parameter $\tau$ on each fine edge $F \in \mathcal{E}_h$ and the coarse scale finite dimensional space $M(F_H)$. For the well-posedness of the scheme, we have the following result: 

\begin{theorem}
Assume that the local spaces satisfy the compatibility conditions of
\cite[Theorem~2.5]{EfendievLazarovShi2015}, and that the stabilization parameter is positive on at least one face
of every fine element and on a nonempty portion of every coarse trace
face.  Then the MsHDG problem has a unique solution.
\end{theorem}

This result is a slightly simplified version of Theorem 2.5 in \cite{EfendievLazarovShi2015}. In this paper, we choose $\tau = 1$ on each $F \in \mathcal{E}_h$, and 
$$W(K) = P_1(K), \boldsymbol{V}(K) = [P_1(K)]^d, M(F_h) = P_1(F_h)$$ for each $K \in \mathcal{T}_h, F_h \in \mathcal{E}_h$. Consequently, the scheme is well-posed for any choice of the coarse space $M(F_H)$.

We briefly recall the role of the coarse trace space in the error
analysis of the standard MsHDG method.  Under the assumptions of
\cite[Theorem~3.8]{EfendievLazarovShi2015} with the above choice of fine-scale local
spaces and with the coarse trace space
\[
    M_H(F_H)=P_\ell(F_H).
\]
If the coarse-block diameter is of order 1, the nonzero
stabilization value is \(\tau=1\), and the single-face stabilization
condition of~\cite{EfendievLazarovShi2015} is satisfied, then Case~1
following Theorem~3.8 gives
\[
    \|u-u_h\|_{L^2(\Omega)}
    \le
    C\left(
       H^{\ell+\frac32}
       \max\{1,h^{1/2}H^{-1}\}
       +h^2
    \right).
\]

An important ingredient in the proof is the approximation of the exact
solution trace by the coarse space in an \(H^{1/2}\)-type norm on the
coarse skeleton.  Schematically, the corresponding contribution takes
the form
\[
    \|u-\mathcal I_Hu\|_{H^{1/2}(\mathcal E_H)},
\]
where \(\mathcal I_Hu\in M_H\) is a suitable coarse interpolant.  This
shows that the approximation quality of \(M_H\), rather than its
polynomial form itself, is the essential issue in the analysis.  The
well-posedness of the MsHDG formulation allows considerable flexibility
in the choice of the coarse trace space, provided that the corresponding
approximation properties are available.  In the present work, however,
we retain the discontinuous piecewise-linear coarse trace spaces and
focus exclusively on the effect of replacing the standard local HDG
solvers by neural-network surrogates.

The stabilization used in the present numerical implementation is
\(\tau=1\) on every fine-grid face.  This choice satisfies the
well-posedness conditions of the MsHDG method and simplifies the
implementation, but it does not satisfy the single-face stabilization
assumption used in the proof of
\cite[Theorem~3.8]{EfendievLazarovShi2015}.  As discussed in that
reference, without the resulting \(H(\operatorname{div})\)-conformity
of the projected flux error, the available estimate contains an
unfavorable factor of \(h^{-1/2}\), including a contribution of order
\(h^{-1/2}H^{\ell+1}\) in the flux analysis.  We therefore do not claim
the optimal global convergence estimate above for the stabilization
used in the present computations.

The objective of the present work is not to establish a new global
convergence theory for this stabilization choice, but to determine
whether neural-network surrogates can replace the repeated fine-scale
local HDG solves required to construct the MsHDG coarse-block
operators.  The standard and neural-network-based methods use the same
fine and coarse discretizations, the same stabilization
\(\tau=1\) on every fine-grid face, and the same global skeleton
coupling.  The only modification is that the local DtN matrices and
source vectors computed by the fine-scale HDG solver are replaced by
their neural-network predictions.  Consequently, comparisons of the
two assembled solutions measure the additional modeling error
introduced by the learned local operators.

The potential benefit of this replacement is computational efficiency.
In the standard MsHDG method, construction of the global skeleton
system requires repeated applications of the inverse of each
fine-scale local HDG matrix to the coarse trace basis functions and the
local source term.  In the neural-network-based method, these local
operators are instead generated by inexpensive forward evaluations of
a trained network.  The numerical experiments in Section~6 therefore
assess both the accuracy of the learned assembly and the resulting
reduction in online assembly and total solution times.

\section{The Upscaling of the MsHDG method}\label{sec:mutiscaleHDG}

In this section, we demonstrate the upscaling structure of the MsHDG scheme and explain the main motivation of the present work.  The global unknown of the method is the trace of the scalar variable on the coarse skeleton.  The fine-scale variables inside each coarse block are eliminated by local HDG solvers, and the global coupling is obtained by enforcing the continuity of the normal numerical flux across the coarse interfaces.

We split the third equation of \eqref{eq:fineHDG-3} by testing separately with functions in $M_h^0$ and $M_H$.  Namely,
\begin{equation}\label{eq:separate}
    \left\langle
        \widehat{\boldsymbol q}_{h,H}\cdot\boldsymbol n,\mu
    \right\rangle_{\partial\mathcal T_h}=0
    \qquad \forall \mu\in M_h^0,
\end{equation}
and
\begin{equation}\label{eq:coarse}
    \sum_i
    \left\langle
        \widehat{\boldsymbol q}_{h,H}\cdot\boldsymbol n,\mu
    \right\rangle_{\partial\Omega_i}=0
    \qquad \forall \mu\in M_H .
\end{equation}
Here $\Omega_i$ denotes a coarse block.  Equation \eqref{eq:separate} is local inside each coarse block, while \eqref{eq:coarse} is the global flux-continuity condition on the coarse skeleton.

For a fixed coarse block $\Omega_i$, define the interior fine skeleton of the block by
\begin{equation}\label{eq:Gamma-i-def}
    \Gamma_i^0:=\partial\mathcal T_h(\Omega_i)\setminus \partial\Omega_i .
\end{equation}
Thus
\begin{equation}\label{eq:boundary-splitting}
    \partial\mathcal T_h(\Omega_i)=\Gamma_i^0\cup \partial\Omega_i .
\end{equation}
To avoid mixing the two trace variables, we use
\[
    \widehat{u}_{h,H}=
    \begin{cases}
        \widehat u_h, & \hbox{on } \Gamma_i^0,\\
        \widehat u_H, & \hbox{on } \partial\Omega_i .
    \end{cases}
\]
The local fine-scale trace $\widehat u_h$ is an unknown on $\Gamma_i^0$, while the coarse trace $\widehat u_H$ is the imposed boundary data on $\partial\Omega_i$.  Given a coarse trace value $\xi_H\in M_H|_{\partial\Omega_i}$, we impose
\begin{equation}\label{eq:local-boundary-xi}
    \widehat u_H=\xi_H
    \qquad \hbox{on } \partial\Omega_i .
\end{equation}
Therefore the hybrid unknown solved by the local problem is only $\widehat u_h$ on $\Gamma_i^0$.

We now replace $\widehat{\boldsymbol{q}}_{h,H}\cdot \boldsymbol n$ with its definition \eqref{eq:trace}. After substituting \eqref{eq:local-boundary-xi} on $\partial\Omega_i$, the second HDG equation may be written in a form that is convenient for the symmetric saddle-point representation below. Indeed, combining the terms involving $\boldsymbol q_h$ over the complete fine-element boundary and integrating by parts elementwise gives
\[
    (\nabla\cdot\boldsymbol q_h,w)_{\mathcal T_h(\Omega_i)}
    +\left\langle \tau(u_h-\widehat u_h),w\right\rangle_{\Gamma_i^0}
    +\left\langle \tau(u_h-\xi_H),w\right\rangle_{\partial\Omega_i}
    =(f,w)_{\mathcal T_h(\Omega_i)}.
\]
Multiplying this identity by $-1$ yields an equivalent equation whose algebraic form is symmetric with the first and third local equations. Thus the local problem on $\Omega_i$ can equivalently be written as follows: find
\[
    (\boldsymbol q_h,u_h,\widehat u_h)
    \in
    \boldsymbol V_h|_{\Omega_i}\times W_h|_{\Omega_i}\times M_h^0|_{\Gamma_i^0}
\]
such that
\begin{subequations}\label{eq:local-expanded-system}
\begin{align}
(\alpha\boldsymbol q_h,\boldsymbol r)_{\mathcal T_h(\Omega_i)}
-(u_h,\nabla\cdot\boldsymbol r)_{\mathcal T_h(\Omega_i)}
+\left\langle
    \widehat u_h,\boldsymbol r\cdot\boldsymbol n
 \right\rangle_{\Gamma_i^0}
&=
-\left\langle
    \xi_H,\boldsymbol r\cdot\boldsymbol n
 \right\rangle_{\partial\Omega_i},
\label{eq:local-expanded-system-a}
\\
-(\nabla\cdot\boldsymbol q_h,w)_{\mathcal T_h(\Omega_i)}
-\left\langle \tau(u_h-\widehat u_h),w\right\rangle_{\Gamma_i^0}
-\left\langle \tau u_h,w\right\rangle_{\partial\Omega_i}
&=
-(f,w)_{\mathcal T_h(\Omega_i)}
-\left\langle \tau\xi_H,w\right\rangle_{\partial\Omega_i},
\label{eq:local-expanded-system-b}
\\
\left\langle
    \boldsymbol q_h\cdot\boldsymbol n+\tau(u_h-\widehat u_h),\mu
\right\rangle_{\Gamma_i^0}
&=0,
\label{eq:local-expanded-system-c}
\end{align}
\end{subequations}
for all
\[
    (\boldsymbol r,w,\mu)
    \in
    \boldsymbol V_h|_{\Omega_i}\times W_h|_{\Omega_i}\times M_h^0|_{\Gamma_i^0}.
\]
In \eqref{eq:local-expanded-system}, all boundary integrals involving the coarse trace have been separated from the integrals over $\Gamma_i^0$ and moved to the right-hand side. The only trace unknown in the local system is therefore $\widehat u_h$ on the interior fine skeleton $\Gamma_i^0$. The sign change in \eqref{eq:local-expanded-system-b} is purely algebraic and does not alter the local HDG solution.

By linearity, the solution of \eqref{eq:local-expanded-system} can be decomposed into the sum of the source-induced part and the trace-induced part,
\begin{equation}\label{eq:split}
    (\boldsymbol q_h,u_h,\widehat u_h)
    =
    (\boldsymbol q_h(f),u_h(f),\widehat u_h(f))
    +
    (\boldsymbol q_h(\xi_H),u_h(\xi_H),\widehat u_h(\xi_H)),
\end{equation}
where the first term is obtained by setting $\xi_H=0$, and the second term is obtained by setting $f=0$. Hence the global system \eqref{eq:coarse} can be written as:
\begin{equation}\label{eq:global_eq}
a_H(\widehat{u}_H, \mu) = l_H(\mu), \quad \forall \mu \in M_H,
\end{equation}
where
\begin{align}\label{eq:global_eq1}
a_H(\widehat{u}_H, \mu) &:=  \sum_i
    \left\langle
        {\boldsymbol q}_{h}(\widehat{u}_H)\cdot \boldsymbol n + \tau (u_h(\widehat{u_H}) - \widehat{u}_H),\mu
    \right\rangle_{\partial\Omega_i}, \\
    l_H(\mu) &:=  - \sum_i
    \left\langle
        {\boldsymbol q}_{h}(f)\cdot \boldsymbol n + \tau u_h(f),\mu
    \right\rangle_{\partial\Omega_i}.
\end{align}

We next write \eqref{eq:local-expanded-system} in matrix form in a block-by-block saddle point form.  Choose bases
\[
    \{\boldsymbol\psi_m\}_{m=1}^{N_q}\subset \boldsymbol V_h|_{\Omega_i},
    \qquad
    \{\varphi_m\}_{m=1}^{N_u}\subset W_h|_{\Omega_i},
    \qquad
    \{\eta_m\}_{m=1}^{N_{\widehat{u}}} \subset M_h^0|_{\Gamma_i^0},
\]
and a basis
\[
    \{\chi_m\}_{m=1}^{N_H}\subset M_H|_{\partial\Omega_i}
\]
for the imposed coarse trace on \(\partial\Omega_i\).  We write the local unknowns as
\[
    \boldsymbol q_h=\sum_m Q_m\boldsymbol\psi_m,
    \qquad
    u_h=\sum_m U_m\varphi_m,
    \qquad
    \widehat u_h=\sum_m \widehat U_m\eta_m,
    \qquad
    \xi_H=\sum_m \Lambda_m\chi_m .
\]
Let
\[
    \boldsymbol Q_i=[Q_m]_{m=1}^{N_q},
    \qquad
    \boldsymbol U_i=[U_m]_{m=1}^{N_u},
    \qquad
    \widehat{\boldsymbol U}_i=[\widehat U_m]_{m=1}^{N_{\widehat{u}}},
    \qquad
    \boldsymbol\lambda_i=[\Lambda_m]_{m=1}^{N_H}.
\]
Then the equivalent form \eqref{eq:local-expanded-system} gives the following symmetric saddle point system:
\begin{equation}\label{eq:local-matrix-system-expanded}
\begin{bmatrix}
    M_i(\alpha) & -B_i^T & G_i^T \\
    -B_i & -\tau T_i & \tau N_i^T \\
    G_i & \tau N_i & -\tau P_i
\end{bmatrix}
\begin{bmatrix}
    \boldsymbol Q_i\\
    \boldsymbol U_i\\
    \widehat{\boldsymbol U}_i
\end{bmatrix}
=
\begin{bmatrix}
    -H_i^T\boldsymbol\lambda_i\\
    -\boldsymbol F_i-\tau J_i\boldsymbol\lambda_i\\
    \boldsymbol 0
\end{bmatrix}.
\end{equation}
Here the blocks are defined by the bilinear forms appearing in \eqref{eq:local-expanded-system}.  More precisely,
\begin{align*}
    [M_i(\alpha)]_{N_q \times N_q}
    &:=(\alpha\boldsymbol\psi_m,\boldsymbol\psi_\ell)_{\mathcal T_h(\Omega_i)},
    &
    [B_i]_{N_u \times N_q}
    &:=(\nabla\cdot\boldsymbol\psi_m,\varphi_\ell)_{\mathcal T_h(\Omega_i)},
    \\
    [G_i]_{N_{\widehat{u}} \times N_q}
    &:=\left\langle \boldsymbol\psi_m\cdot\boldsymbol n,\eta_\ell\right\rangle_{\Gamma_i^0},
    &
    [T_i]_{N_u \times N_u}
    &:=\left\langle \varphi_m,\varphi_\ell\right\rangle_{\partial\mathcal T_h(\Omega_i)},
    \\
    [N_i]_{N_{\widehat{u}} \times N_u}
    &:=\left\langle \varphi_\ell, \eta_m\right\rangle_{\Gamma_i^0},
    &
    [P_i]_{N_{\widehat{u}} \times N_{\widehat{u}}}
    &:=\left\langle \eta_m,\eta_\ell\right\rangle_{\Gamma_i^0},
    \\
    [H_i]_{ N_H \times N_q}
    &:=\left\langle \boldsymbol\psi_m\cdot\boldsymbol n,\chi_\ell\right\rangle_{\partial\Omega_i},
    &
    [J_i]_{N_u \times N_H}
    &:=\left\langle \chi_m,\varphi_\ell\right\rangle_{\partial\Omega_i}.
\end{align*}
The vector \(\boldsymbol F_i\) is defined by
\[
    [\boldsymbol F_i]_{l=1}^{N_u}=(f,\varphi_\ell)_{\mathcal T_h(\Omega_i)}.
\]
In this notation the first block row corresponds to \eqref{eq:local-expanded-system-a}, the second block row corresponds to the integrated-by-parts form \eqref{eq:local-expanded-system-b}, and the third block row corresponds to the interior flux-continuity constraint \eqref{eq:local-expanded-system-c}. Notice that the imposed coarse trace \(\xi_H\) only appears on the right-hand side through the two boundary contributions \(-H_i^T\boldsymbol\lambda_i\) and \(-\tau J_i\boldsymbol\lambda_i\).

For compactness, we denote the saddle point matrix in \eqref{eq:local-matrix-system-expanded} by
\begin{equation}\label{eq:local-saddle-matrix-definition}
    A_i(\alpha)
    :=
    \begin{bmatrix}
        M_i(\alpha) & -B_i^T & G_i^T \\
        -B_i & -\tau T_i & \tau N^T_i \\
        G_i & \tau N_i & - \tau P_i
    \end{bmatrix}_{\mathcal{N}_h \times \mathcal{N}_h},
\,  \,
    L_i
    :=
    \begin{bmatrix}
        -H_i^T\\
        -\tau J_i\\
        \bf{0}
    \end{bmatrix}_{\mathcal{N}_h \times N_H},
    \, \,
    \boldsymbol b_i(f)
    :=
    \begin{bmatrix}
        \bf{0}\\
        -\boldsymbol F_i\\
        \bf{0}
    \end{bmatrix}_{\mathcal{N}_h \times 1}.
\end{equation}
Since \(M_i(\alpha)\), \(T_i\), and \(P_i\) are symmetric and the off-diagonal blocks occur in transpose pairs, \(A_i(\alpha)\) is a symmetric saddle point matrix. Its invertibility follows from the well-posedness of the corresponding local HDG problem.
Here $\mathcal{N}_h = N_u + N_q + N_{\widehat{u}}$ which is the number of DoFs of the fine-scale system on each coarse subdomain $\Omega_i$.
Then \eqref{eq:local-matrix-system-expanded} can be written as
\begin{equation}\label{eq:local-matrix-system-compact}
    A_i(\alpha)
    \begin{bmatrix}
        \boldsymbol Q_i\\
        \boldsymbol U_i\\
        \widehat{\boldsymbol U}_i
    \end{bmatrix}
    =
    \boldsymbol b_i(f)+L_i\boldsymbol\lambda_i .
\end{equation}
Therefore the local solver is
\begin{equation}\label{eq:local-solver-matrix-expanded}
    \begin{bmatrix}
        \boldsymbol Q_i\\
        \boldsymbol U_i\\
        \widehat{\boldsymbol U}_i
    \end{bmatrix}
    =
    \underbrace{A_i(\alpha)^{-1}\boldsymbol b_i(f)}_{(\boldsymbol{q}_h(f), u_h(f), \widehat{u}_h(f))}
    +\underbrace{A_i(\alpha)^{-1}L_i\boldsymbol\lambda_i}_{(\boldsymbol{q}_h(\xi_H), u_h(\xi_H), \widehat{u}_h(\xi_H))} .
\end{equation}
Thus, the local upscaling step requires applying the inverse of the fine-scale saddle point matrix \(A_i(\alpha)\).

The quantity needed for the global assembly in \eqref{eq:coarse} is the numerical flux $\boldsymbol{q}_{h} \cdot \boldsymbol{n} + \tau (u_h - \xi_H)$ on the coarse boundary $\partial\Omega_i$.  
Let the matrix $[C_i]_{N_H \times \mathcal{N}_h} = [H_i \quad \tau J^T_i \quad  \bf{0}]$ denote the restriction mapping of $(\boldsymbol{q}_h, u_h)|_{\Omega_i} \rightarrow \boldsymbol{q}_h \cdot \boldsymbol{n} + \tau u_h|_{\partial \Omega_i}$ and $[D_i]_{N_H \times N_H} = \langle \chi_m, \chi_l \rangle_{\partial \Omega_i}$. With the symmetric convention \eqref{eq:local-saddle-matrix-definition}, we have
\[
    C_i=-L_i^T.
\]
Then the assembly matrices and the right-hand-side vectors for the global system \eqref{eq:global_eq}, \eqref{eq:global_eq1} for each coarse block $\Omega_i$ can be written as:
\begin{subequations}\label{eq:local-flux-vector-pre}
\begin{align}
   S_i(\alpha) &:=  [a_H(\chi_m, \chi_l)_{\Omega_i}]_{N_H \times N_H} 
    =
    C_iA_i(\alpha)^{-1}L_i - \tau D_i
    =-L_i^TA_i(\alpha)^{-1}L_i-\tau D_i, \\
    g_i(\alpha,f)
     & :=
    C_iA_i(\alpha)^{-1}\boldsymbol{b}_i(f).
\end{align}
\end{subequations}

The matrix $S_i(\alpha)$ is the discrete Dirichlet-to-Neumann matrix on the coarse block $\Omega_i$. It maps the imposed coarse trace degrees of freedom on $\partial\Omega_i$ to the corresponding normal numerical flux degrees of freedom on $\partial\Omega_i$ with $f = 0$. Since $A_i(\alpha)^{-1}$ and $D_i$ are symmetric, the representation above shows directly that $S_i(\alpha)$ is symmetric. The vector $\boldsymbol g_i(\alpha,f)$ is the source-induced flux contribution.

We can now express the global coarse system \eqref{eq:coarse} in the matrix form.  Let
\[
    M_H=\operatorname{span}\{\chi_1,\ldots,\chi_\mathcal{M}\},
    \qquad
    \widehat{u}_H=\sum_{j=1}^\mathcal{M} \widehat{U}_j \chi_j,
    \qquad
    \widehat{\boldsymbol{U}}_H=(\widehat{U}_1,\ldots,\widehat{U}_{\mathcal{M}})^T .
\]
For each coarse block $\Omega_i$, let $R_i$ be the restriction matrix from the global trace vector $\widehat{\boldsymbol U}_H$ to the local trace vector on $\partial\Omega_i$:
\begin{equation}\label{eq:local-restriction-expanded}
    \boldsymbol\lambda_i=R_i\widehat{\boldsymbol U}_H .
\end{equation}
By inserting \eqref{eq:local-flux-vector-pre} into the global skeleton equations \eqref{eq:global_eq}, we obtain the global system in the matrix form as: 
\begin{equation}
    \left(
      \sum_i R_i^TS_i(\alpha)R_i
    \right)\widehat{\boldsymbol U}_H
    =
    -\sum_i R_i^Tg_i(\alpha,f).
\end{equation}

This is the matrix form of the MsHDG upscaling procedure.  Although the final system only contains the coarse trace unknowns, its exact assembly requires constructing the local DtN matrix $S_i(\alpha)$ on every coarse block.

The computational bottleneck is apparent from \eqref{eq:local-flux-vector-pre}.  To compute $S_i(\alpha)$ exactly, one must apply the inverse of the fine-scale local matrix $A_i(\alpha)$ to all trace basis inputs on $\partial\Omega_i$.  The matrix $A_i(\alpha)$ is assembled from the heterogeneous coefficient $\alpha=\kappa^{-1}$ on the fine mesh of $\Omega_i$.  Consequently, the DtN matrix $S_i(\alpha)$ depends on the local permeability field in a highly nonlinear way through the inverse $A_i(\alpha)^{-1}$.  Therefore, for a heterogeneous medium, this fine-scale matrix inversion has to be performed separately on every coarse block.  This local construction of the DtN matrices becomes the major computational bottleneck of the MsHDG scheme in practice.


The main idea of the present work is to replace the repeated
construction of the local coarse-block operators by a neural-network
surrogate.  In this study, the geometry and size of the coarse blocks
are fixed, and the source term is prescribed as \(f\equiv1\).
Consequently, neither the coarse-block geometry nor the source term is
treated as a neural-network input.  The only varying input is the local
permeability field \(\kappa|_{\Omega_i}\), or equivalently
\(\alpha|_{\Omega_i}=\kappa^{-1}|_{\Omega_i}\).  We therefore consider
the coefficient-to-operator map
\[
    \kappa|_{\Omega_i}
    \longmapsto
    \bigl(
        S_i(\kappa^{-1}),
        g_i(\kappa^{-1},1)
    \bigr).
    \tag{4.18}
\]
Although the Darcy problem is linear with respect to the trace and
source data, the local operators depend nonlinearly on the permeability
field through the inverse of the heterogeneous fine-scale matrix
\(A_i(\kappa^{-1})\).

We train a neural network \(\mathcal N_\theta\) to approximate this map,
\[
    \bigl(
        S_i(\kappa^{-1}),
        g_i(\kappa^{-1},1)
    \bigr)
    \approx
    \mathcal N_\theta\!\left(\kappa|_{\Omega_i}\right).
    \tag{4.19}
\]
Since all coarse blocks used in the numerical experiments are
geometrically congruent, their coefficient fields and local operators
can be represented on a common reference block.  After training, the
local DtN matrix and source-induced flux vector are obtained by a
forward evaluation of the network rather than by repeated applications
of \(A_i(\kappa^{-1})^{-1}\).  The predicted operators are then
assembled using the same global MsHDG skeleton equations as in the
standard method.

\section{Neural-network-based local solvers}\label{sec:NN solvers}

We now describe the reference-block discretization, training-data
generation, and neural-network architecture used to approximate the
local operators introduced in Section \ref{sec:mutiscaleHDG}. For simplicity, throughout the numerical experiments we fix the source
term as $f \equiv 1$
and the coarse partition is a uniform partition of $\Omega = [0,1]^2$ with coarse squares.  Consequently, the only
varying input parameter of the neural network is the local heterogeneous
coefficient
\[
    \alpha|_{\Omega_i}=\kappa^{-1}|_{\Omega_i}.
\]
The network is trained to approximate the coefficient-to-local-operator
map
\[
    \alpha|_{\Omega_i}
    \longmapsto
    \bigl(S_i(\alpha),\,g_i(\alpha,1)\bigr).
\]
Here we choose \(S_i(\alpha)\) and \(g_i(\alpha,1)\), rather than the full
inverse matrix \(A_i(\alpha)^{-1}\), as the neural-network outputs.
Predicting the full inverse would therefore
require \(\mathcal{N}_h^2\) output entries, whereas predicting the symmetric matrix
\(S_i(\alpha)\) together with \(g_i(\alpha,1)\in\mathbb{R}^{N_H}\)
requires only $\frac{N_H(N_H+1)}{2}+N_H$ independent output components.  Thus, learning the reduced quantities
\(\bigl(S_i(\alpha),g_i(\alpha,1)\bigr)\) substantially decreases the
output dimension while providing all information needed for the global
assembly.

A central advantage of the proposed approach is that the neural network can be trained on a reference coarse block to approximate the local solution operator, or more specifically the Dirichlet-to-Neumann (DtN) map associated with that block.  Once trained, the same reference-block model can be applied to other coarse blocks after the appropriate geometric and coefficient representation has been supplied.

\subsection{Reference-block discretization and data generation}

All two-dimensional training data are generated on the reference coarse
block $\widehat K=[0,1]^2.$
We consider a two-phase heterogeneous medium formed by materials with
permeabilities \(\kappa_{\mathrm{low}}\) and
\(\kappa_{\mathrm{high}}\).  For simplicity, we normalize the
low-permeability phase by setting $ \kappa_{\mathrm{low}}=1.$
The numerical experiments consider the two choices
\[
    \kappa_{\mathrm{high}}=10
    \qquad\text{and}\qquad
    \kappa_{\mathrm{high}}=10^4,
\]
which will be referred to as the low-contrast and high-contrast regimes,
respectively.

\begin{figure}[H]
    \centering
    \begin{subfigure}[t]{0.48\linewidth}
        \centering
        \includegraphics[width=\linewidth]
        {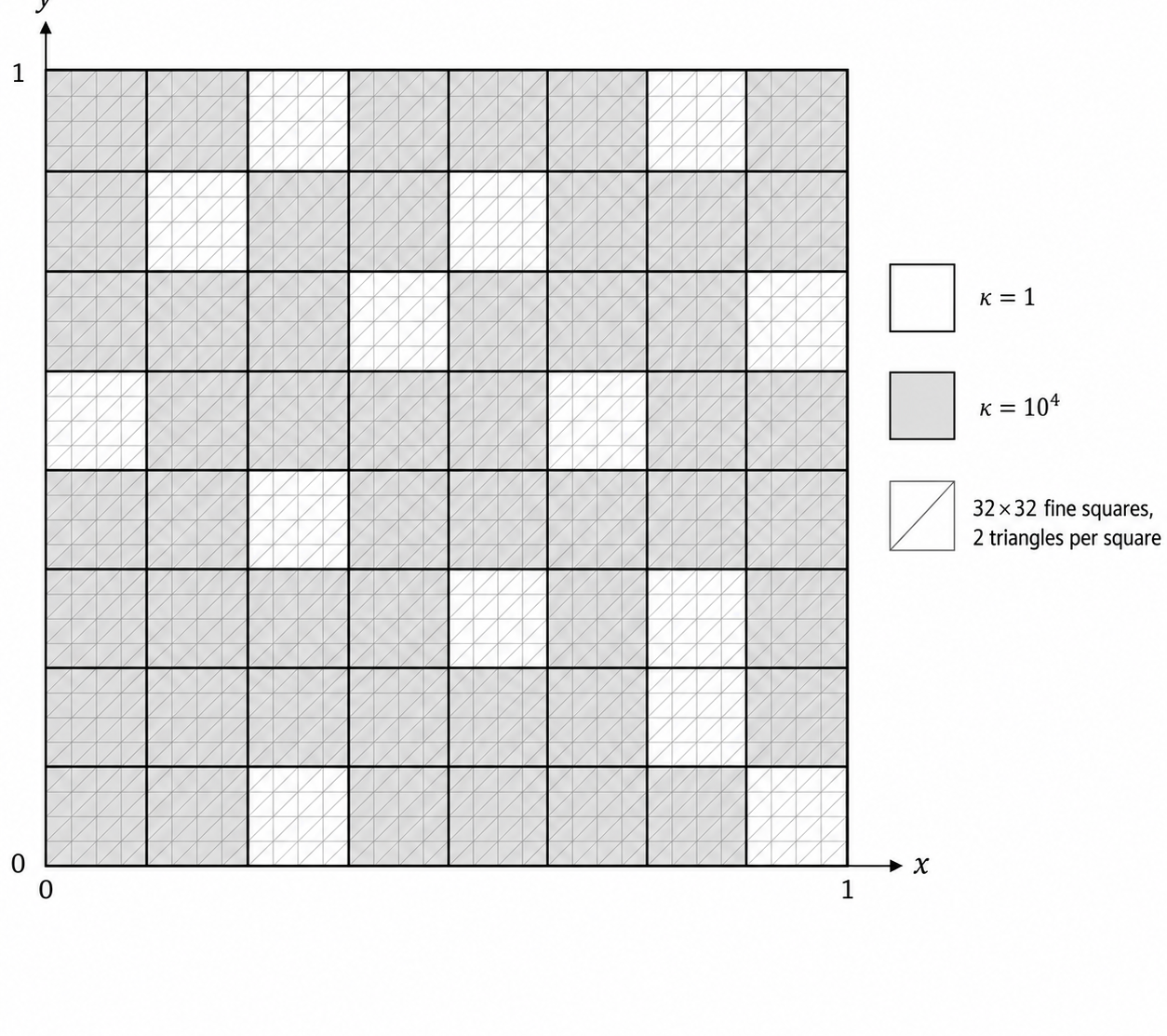}
        \caption{Example A.}
    \end{subfigure}
    \hfill
    \begin{subfigure}[t]{0.48\linewidth}
        \centering
        \includegraphics[width=\linewidth]
        {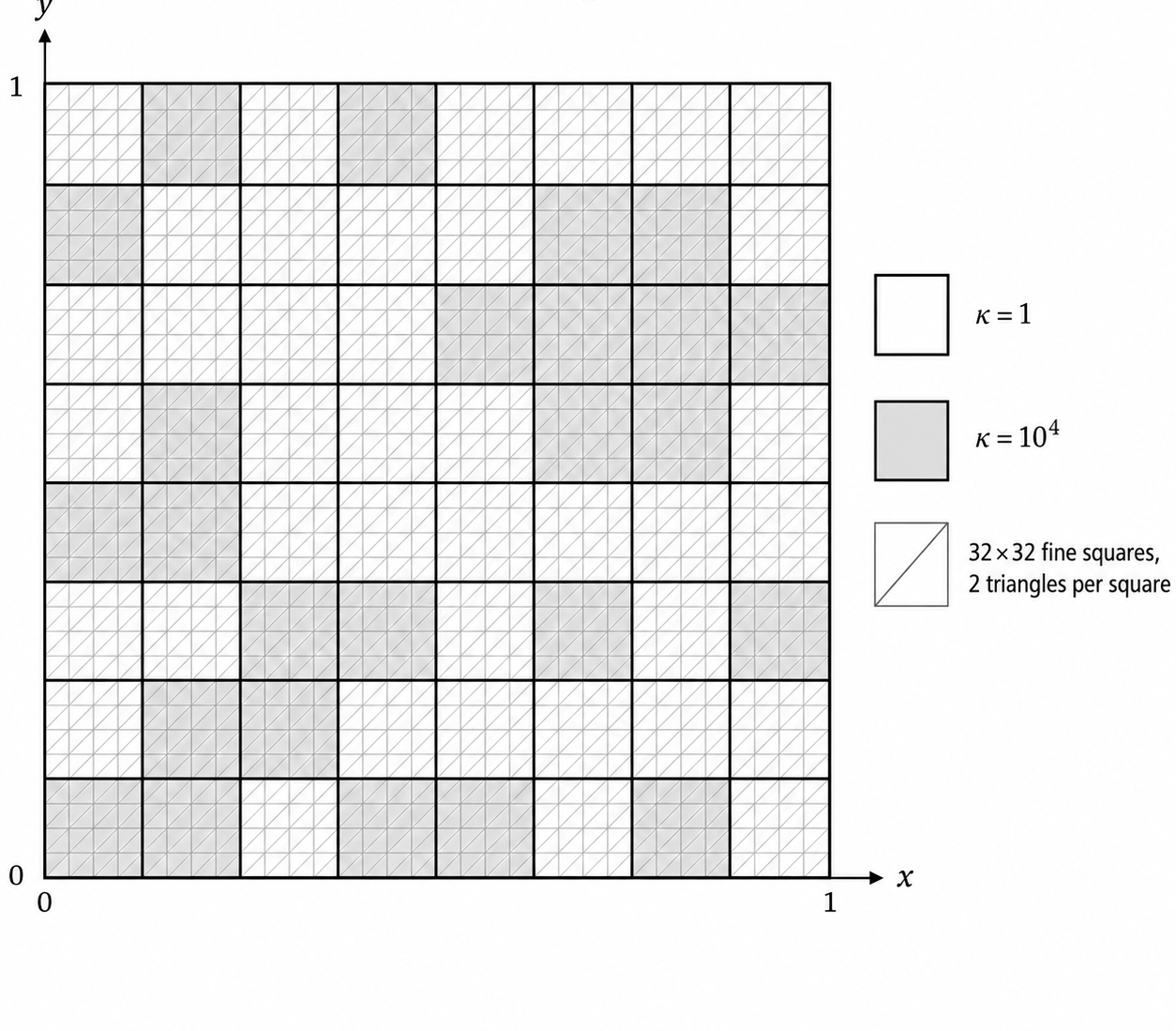}
        \caption{Example B.}
    \end{subfigure}
    \caption{Examples of random \(8\times8\) reference-block coefficient
    distributions overlaid with the fine triangular mesh.  The white and
    gray cells represent the phases
    \(\kappa_{\mathrm{low}}\) and \(\kappa_{\mathrm{high}}\),
    respectively.}
    \label{fig:kappa-examples}
\end{figure}

The coefficient field is represented on a uniform \(8\times8\) feature
grid on \(\widehat K\).  For each coefficient realization, the 64 feature
cells are sampled independently from a Bernoulli two-phase distribution,
\[
    \mathbb P\!\left(\kappa=\kappa_{\mathrm{low}}\right)=\frac12,
    \qquad
    \mathbb P\!\left(\kappa=\kappa_{\mathrm{high}}\right)=\frac12.
\]
Thus, each coefficient realization is described by an independently sampled
binary \(8\times8\) material pattern. The two
examples in Figure~\ref{fig:kappa-examples} illustrate typical coefficient
fields used in the generation of the local training and validation datasets.
The same reference geometry and Bernoulli sampling law are used for both
contrast regimes.  


The fine mesh on \(\widehat K\) is obtained by partitioning the reference
block into \(32\times32\) uniform squares and then dividing each square
into two triangles.  Consequently, each coefficient feature cell
contains a \(4\times4\) patch of fine squares, or equivalently \(32\)
fine triangles.  The neural network therefore receives a fixed
\(8\times8\) coefficient image, while the corresponding local operators
are computed by HDG solves on the finer triangular mesh.

We next define the coarse trace spaces on
\(\partial\widehat K\).  For each coarse edge
\(\widehat F\subset\partial\widehat K\) and each trace level
\(n\geq0\), let
\[
    \mathcal E_H^{(n)}(\widehat F)
      =
      \{\widehat F_j^{(n)}\}_{j=1}^{2^n}
\]
be the uniform partition of \(\widehat F\) into \(2^n\) subintervals.
The corresponding trace space is
\[
    M_H^{(n)}(\widehat F)
      :=
      \left\{
        \mu\in L^2(\widehat F):
        \mu|_{\widehat F_j^{(n)}}
        \in P_1(\widehat F_j^{(n)}),
        \quad j=1,\ldots,2^n
      \right\}.
\]
No continuity is imposed between adjacent subintervals, since every
subinterval contributes two discontinuous piecewise-linear degrees of
freedom, $\dim M_H^{(n)}(\widehat F)=2^{n+1}.$
The trace space on the boundary of the reference block is
\[
    M_H^{(n)}(\partial\widehat K)
      =
      \bigoplus_{\widehat F\subset\partial\widehat K}
      M_H^{(n)}(\widehat F).
\]
The total dimension of the trace input to the local DtN map is therefore
\[
    M(n)
      :=
      \dim M_H^{(n)}(\partial\widehat K)
      =
      4\cdot 2^{n+1}
      =
      2^{n+3}.
\]
Thus, increasing \(n\) uniformly refines the trace partition on each
coarse edge and enriches the coarse trace space.

On the reference block \(\widehat K\), for each sampled coefficient field
\(\alpha=\kappa^{-1}\), the corresponding training labels are generated
by the MsHDG local solver \eqref{eq:local-flux-vector-pre}, with the
source term fixed as \(f\equiv1\).  This produces the local
coefficient-to-operator map
\[
    \alpha|_{\widehat K}
    \longmapsto
    \bigl(S_{\widehat K}(\alpha),
          g_{\widehat K}(\alpha,1)\bigr),
\]
where \(S_{\widehat K}(\alpha)\) and
\(g_{\widehat K}(\alpha,1)\) are the local contributions used in the
assembly of the global coarse-scale system.

The permeability input is normalized cellwise according to
\[
    X
      =
      \frac{\log_{10}\kappa}
           {\log_{10}\kappa_{\mathrm{high}}}.
\]
For either of the two binary contrast regimes considered here, this
normalization maps the two material phases to
\[
    \kappa=1 \longmapsto X=0,
    \qquad
    \kappa=\kappa_{\mathrm{high}} \longmapsto X=1.
\]
Hence, each coefficient realization is represented by a one-channel
\(8\times8\) normalized permeability image.  Equivalently, the
neural-network input is a vector in \(\mathbb R^{64}\), and its
dimension is independent of the trace level \(n\).

For a fixed trace level \(n\), the coarse trace space on
\(\partial\widehat K\) has dimension
\[
    M(n)
      =
      \dim M_H^{(n)}(\partial\widehat K)
      =
      2^{n+3}.
\]
Accordingly,
\[
    S(\alpha)\in\mathbb R^{M(n)\times M(n)},
    \qquad
    g(\alpha,1)\in\mathbb R^{M(n)}.
\]
By the symmetric representation in \eqref{eq:local-flux-vector-pre},
\(S(\alpha)\) is symmetric, so only its upper-triangular entries are retained.  The vector \(g(\alpha,1)\) is stored in full.  Therefore, for
each coefficient sample, the neural-network output is the concatenated
vector
\[
    Y^{(n)}
      =
      \left(
        \operatorname{upper}\bigl(S(\alpha)\bigr),
        g(\alpha,1)
      \right),
\]
with total dimension
\[
    N_{\mathrm{out}}(n)
      =
      \frac{M(n)\bigl(M(n)+1\bigr)}{2}
      +M(n).
\]

The neural-network input dimension, coarse trace dimension, and output
dimension for the trace levels used in the numerical experiments are
summarized in Table~\ref{tab:trace-dims}.  The input dimension remains
fixed at \(64\), whereas the trace-space and output dimensions increase
with the trace enrichment level \(n\).

\begin{table}[H]
    \centering
    \caption{Neural-network input, coarse trace, and output dimensions
    used in the two-dimensional experiments.}
    \label{tab:trace-dims}
    \begin{tabular}{cccc}
        \toprule
        Trace level \(n\)
        & NN input dimension
        & Coarse trace dimension \(M(n)\)
        & NN output dimension \(N_{\mathrm{out}}(n)\)
        \\
        \midrule
        0 & 64 & 8  & 44  \\
        1 & 64 & 16 & 152 \\
        2 & 64 & 32 & 560 \\
        \bottomrule
    \end{tabular}
\end{table}

For each contrast regime, we generate
$N_{\mathrm{samples}}=5000$
independent local reference-block coefficient realizations. The data are randomly divided into an $80\%$ training set and a $20\%$ validation set, with the latter used for checkpoint selection and training diagnostics. Independent global physical-domain realizations, evaluated only after training, are used for the numerical tests reported in Section~6. The same coefficient representation, trace spaces, output ordering, and fine-scale HDG discretization are used in both contrast regimes.

All physical coarse blocks used in the numerical experiments are congruent
squares obtained from the reference block \(\widehat K=[0,1]^2\) by an
affine scaling and translation. The local coefficient field on each
physical block is therefore represented in the same reference coordinates
used for training. The corresponding local DtN matrix and source vector
are transferred between the reference and physical blocks using the
deterministic geometric scaling associated with this affine map. Since the
coarse-block geometry is fixed throughout the present experiments, no
geometric parameters are included in the neural-network input, and the
same reference-block network can be applied to every coarse block.

\subsection{Network architecture and training}

For each contrast regime and each trace level $n=0,1,2$, we train a
separate neural network using the datasets described in the preceding
subsection. The neural-network input and output dimensions are those
listed in Table~1.

Since the permeability input is represented as an $8\times 8$ material
image, we use a convolutional neural network to extract spatial features
of the coefficient field. The network consists of three convolutional
layers followed by two fully connected hidden layers. Each convolution
uses a $3\times 3$ kernel with unit padding, so that the spatial
resolution remains $8\times 8$ throughout the convolutional part of the
network. We use the Gaussian error linear unit (GELU) activation
function~\cite{HendrycksGimpel2016}, and the final layer is linear since
the entries of $S(\alpha)$ and $g(\alpha,1)$ are real-valued. The
architecture is summarized in Table~\ref{tab:NNarchitecture}.

\begin{table}[ht]
\centering
\caption{Neural-network architecture used in the two-dimensional
contrast experiments.}
\label{tab:NNarchitecture}
\begin{tabular}{ll}
\hline
Stage & Layer description\\
\hline
Input & $1\times 8\times 8$ normalized coefficient image\\
Convolution 1 & $3\times3$ convolution, $1\to32$ channels, padding $1$, GELU\\
Convolution 2 & $3\times3$ convolution, $32\to64$ channels, padding $1$, GELU\\
Convolution 3 & $3\times3$ convolution, $64\to64$ channels, padding $1$, GELU\\
Flatten & $64\times8\times8=4096$ features\\
Fully connected 1 & linear $4096\to512$, GELU\\
Fully connected 2 & linear $512\to512$, GELU\\
Output & linear $512\to N_{\mathrm{out}}(n)$\\
\hline
\end{tabular}
\end{table}

The principal training objective is to reproduce the local DtN matrix
and source vector. Let
\[
    \bigl(\widehat S_\theta(\alpha),
          \widehat g_\theta(\alpha,1)\bigr)
\]
denote the quantities reconstructed from the network output. We first
use the normalized data loss
\begin{equation}
\mathcal L_{\mathrm{data}}
=
\frac{\|\widehat S_\theta(\alpha)-S(\alpha)\|_F^2}
     {\|S(\alpha)\|_F^2+\varepsilon}
+
\frac{\|\widehat g_\theta(\alpha,1)-g(\alpha,1)\|_2^2}
     {\|g(\alpha,1)\|_2^2+\varepsilon},
\label{eq:loss-data}
\end{equation}
where $\varepsilon=10^{-12}$ prevents division by very small
denominators.

Since the DtN matrix enters the global MsHDG system primarily through
its action on coarse trace functions, the entrywise loss is supplemented
by structure-aware terms. For a collection of representative trace
vectors $V\subset\mathbb R^{M(n)}$, we define the relative operator-action
loss
\begin{equation}
\mathcal L_{\mathrm{act}}(V)
=
\frac{1}{|V|}
\sum_{v\in V}
\frac{\|(\widehat S_\theta-S)v\|_2^2}
     {\|Sv\|_2^2+\varepsilon},
\label{eq:loss-action}
\end{equation}
and the corresponding energy-consistency loss
\begin{equation}
\mathcal L_{\mathrm{en}}(V)
=
\frac{1}{|V|}
\sum_{v\in V}
\frac{\bigl|v^T(\widehat S_\theta-S)v\bigr|^2}
     {|v^TSv|^2+\varepsilon}.
\label{eq:loss-energy}
\end{equation}
Finally, the homogeneous DtN operator annihilates the constant trace.
We therefore include the nullspace penalty
\begin{equation}
\mathcal L_{\mathrm{null}}
=
\|\widehat S_\theta \mathbf 1\|_2^2,
\label{eq:loss-null}
\end{equation}
where $\mathbf 1$ denotes the coefficient vector representing the
constant trace on $\partial\widehat K$.

The probe vectors are projected onto the constant-free trace subspace
and normalized before the losses are evaluated. We use three probe
families. The smooth family consists of low-order side modes,
side-average and adjacent-side combinations, together with $32$ random
low-frequency combinations. This construction gives $54$, $78$, and
$126$ smooth probes for $n=0,1,2$, respectively. We also use four fully
random constant-free traces and $16$ solution-representative probes
formed from low-order polynomial, gradient, curvature, and correlated
side patterns. All probes are local to the reference block; no exact
global solution or global skeleton vector is used during training.

The action terms associated with the smooth, random, and
solution-representative probe families are assigned weights $0.45$,
$0.008$, and $0.10$, respectively. The corresponding energy penalties
are given substantially smaller weights, of order $10^{-5}$--$10^{-4}$,
and the nullspace penalty has weight $0.01$. Thus the normalized data
loss remains the dominant term, while the auxiliary losses encourage
the predicted DtN matrix to reproduce its action on trace directions
that are relevant to the subsequent global assembly.

The network parameters are optimized using AdamW~\cite{LoshchilovHutter2019}
with an initial learning rate of $3\times10^{-4}$ and weight decay
$10^{-4}$. The learning rate is reduced using a cosine-annealing
schedule~\cite{LoshchilovHutter2017}. We use a batch size of $128$ and
random seed $123$. The maximum training horizons are $600$, $900$, and
$1200$ epochs for trace levels $n=0,1,2$, respectively, and the
validation subset is used for checkpoint selection.

The same network architecture, output ordering, loss construction, and
optimization parameters are used for the two contrast regimes
\[
    \kappa\in\{1,10\}
    \qquad\text{and}\qquad
    \kappa\in\{1,10^4\}.
\]
Separate networks are trained for the two regimes because the
corresponding local HDG operators depend strongly on the value of
$\kappa_{\mathrm{high}}$. This common experimental setup allows the
effect of coefficient contrast on the learned local operators to be
compared directly.

\section{Numerical experiments}
\label{sec:numerical-experiments}

The purpose of this section is to assess the accuracy and online computational efficiency of the proposed neural-network-based multiscale HDG method.  Throughout the numerical experiments, the source term is fixed as $f\equiv1$, and homogeneous Dirichlet boundary conditions are imposed on $\partial\Omega$.

We distinguish three discrete solutions.  Let $u_h^{\mathrm{REF}}$ denote the standard HDG solution computed on the complete fine mesh, let $u_H^{\mathrm{MS}}$ denote the standard MsHDG solution assembled from the local operators computed by the fine-scale HDG solver, and let $u_H^{\mathrm{NN}}$ denote the NN-MsHDG solution assembled from the neural-network-predicted local operators.  The subscript $H$ emphasizes that the globally coupled unknowns of the two multiscale methods belong to the coarse trace space, although their interior solutions are reconstructed on the fine mesh.

Unless stated otherwise, all solution errors are measured in the $L^2(\Omega)$ norm.  We use the following three relative errors:
\begin{equation}
    E_{\mathrm{MS}}
    =
    \frac{\|u_h^{\mathrm{REF}}-u_H^{\mathrm{MS}}\|_{L^2(\Omega)}}
         {\|u_h^{\mathrm{REF}}\|_{L^2(\Omega)}}.
    \label{eq:MS-reference-error}
\end{equation}
\begin{equation}
    E_{\mathrm{NN}}
    =
    \frac{\|u_h^{\mathrm{REF}}-u_H^{\mathrm{NN}}\|_{L^2(\Omega)}}
         {\|u_h^{\mathrm{REF}}\|_{L^2(\Omega)}}.
    \label{eq:NN-reference-error}
\end{equation}
\begin{equation}
    E_{\mathrm{model}}
    =
    \frac{\|u_H^{\mathrm{NN}}-u_H^{\mathrm{MS}}\|_{L^2(\Omega)}}
         {\|u_H^{\mathrm{MS}}\|_{L^2(\Omega)}}.
    \label{eq:NN-modeling-error}
\end{equation}
The error $E_{\mathrm{MS}}$ measures the approximation error of the underlying MsHDG discretization, while $E_{\mathrm{NN}}$ measures the total error of the accelerated method.  Their comparison shows whether NN-MsHDG retains the accuracy of the standard MsHDG method.  The quantity $E_{\mathrm{model}}$ isolates the additional error introduced by replacing the local HDG solves with the neural-network surrogate.

The online computational cost is measured after the offline data generation and training stages have been completed.  We distinguish the time $T_{\mathrm{asm}}$ required to generate all local coarse-block operators and assemble the global skeleton matrix and right-hand side from the total online time $T_{\mathrm{online}}$, which additionally includes the global trace solve and the local reconstruction of the physical solution.  The corresponding speedup factors are
\begin{equation}
    \mathcal S_{\mathrm{asm}}
    =
    \frac{T_{\mathrm{asm}}^{\mathrm{MS}}}
         {T_{\mathrm{asm}}^{\mathrm{NN}}},
    \qquad
    \mathcal S_{\mathrm{online}}
    =
    \frac{T_{\mathrm{online}}^{\mathrm{MS}}}
         {T_{\mathrm{online}}^{\mathrm{NN}}}.
    \label{eq:timing-speedups}
\end{equation}


\subsection{Computational environment and timing protocol}
\label{subsec:computational-environment}

All computations were performed on a desktop workstation equipped with an Intel Core i7-14700F processor, 32 GB RAM, and an NVIDIA GeForce RTX 4060 GPU with 8 GB memory. The HDG and MsHDG solvers were implemented in MATLAB R2026a, while neural-network training and inference used Python 3.11.7 and PyTorch 2.2.2. MATLAB computations were performed in double precision, whereas neural-network computations used single precision.

The reported two-dimensional timing results were obtained on the same workstation using 15 repeated runs for each trace level and contrast. MATLAB parallel processing was disabled during these tests, although its default internal multithreading was retained. Offline dataset generation and neural-network training were excluded from the reported wall-clock times, so the comparison reflects the online computational cost of standard MsHDG and NN-MsHDG after training.


We first use a one-dimensional problem to verify the local learning and global assembly procedure in a setting where the homogeneous DtN map is available analytically.  We then present two-dimensional experiments for the lower-contrast regime $\kappa\in\{1,10\}$ and the high-contrast regime $\kappa\in\{1,10^4\}$.

\subsection{One-dimensional verification}

We first consider a one-dimensional experiment as a simple verification
of the learning and assembly procedure. The main advantage of this setting
is that the homogeneous local Dirichlet-to-Neumann (DtN) operator can be
written explicitly, so that the quantity learned by the neural network has
a direct analytical interpretation.

Let \(K=(x_L,x_R)\) be a coarse interval of length \(H_K\). The homogeneous
local problem satisfies
\[
    \alpha q + u' = 0,
    \qquad
    q' = 0,
\]
where \(\alpha=\kappa^{-1}\). Since \(q\) is constant on \(K\), integration
gives
\[
    u(x_R)-u(x_L)
    =
    -q\int_K \kappa(x)^{-1}\,dx.
\]
Hence,
\[
    q
    =
    -c_K\bigl(u(x_R)-u(x_L)\bigr),
\]
where
\begin{equation}
    c_K
    =
    \left(
        \int_K \kappa(x)^{-1}\,dx
    \right)^{-1}
    \label{eq:one-dimensional-conductance}
\end{equation}
is the effective conductance of the interval. Up to the normal-flux sign
convention, the corresponding local DtN matrix is
\begin{equation}
    S_K
    =
    c_K
    \begin{bmatrix}
        1 & -1\\
        -1 & 1
    \end{bmatrix}.
    \label{eq:one-dimensional-DtN}
\end{equation}

In the numerical experiment, each coarse interval is divided into $N_k$
uniform coefficient cells. If \(\kappa_i\) denotes the permeability on
the \(i\)-th cell, then
\[
    \int_K \kappa(x)^{-1}\,dx
    =
    \frac{H_K}{N_k}
    \sum_{i=1}^{N_k} \frac{1}{\kappa_i},
\]
and therefore
\[
    c_K
    =
    \frac{\kappa_{\mathrm{eff},K}}{H_K},
    \qquad
    \kappa_{\mathrm{eff},K}
    =
    \left(
        \frac{1}{N_k}
        \sum_{i=1}^{N_k} \frac{1}{\kappa_i}
    \right)^{-1}.
\]
Thus, \(\kappa_{\mathrm{eff},K}\) is the harmonic mean of the $N_k$
cellwise permeability values. On the unit reference interval,
\(H_K=1\), the effective conductance coincides with
\(\kappa_{\mathrm{eff},K}\). Consequently, the neural network is trained
on the reference interval to predict the effective permeability, while
the conductance on a physical coarse interval is recovered through the
deterministic scaling \(H_K^{-1}\).

Moreover, because \(\kappa\) is piecewise constant on coefficient-aligned
cells, the homogeneous solution is piecewise linear. These functions
belong to the chosen \(P_1\) HDG spaces, so the fine-scale HDG local solver
reproduces the analytical homogeneous DtN operator. The one-dimensional
problem therefore provides a controlled test in which the learned scalar
quantity can be compared directly with its exact value.

The cellwise permeability values are sampled independently from
\[
    \kappa_i\in\{1,10\}
    \qquad\text{or}\qquad
    \kappa_i\in\{1,10^4\},
    \qquad i=1,\ldots,N_\kappa.
\] For a sample \(s\), the neural-network input is the logarithmic coefficient
vector
\[
    X^{(s)}
    =
    \bigl[
        \log \kappa_1^{(s)},
        \ldots,
        \log \kappa_{N_k}^{(s)}
    \bigr],
\]
and the scalar learning target is
\[
    Y^{(s)}
    =
    \log \kappa_{\mathrm{eff}}^{(s)}.
\]
After inference, the predicted effective permeability is recovered by
\[
    \kappa_{\mathrm{eff}}^{\mathrm{NN}}
    =
    \exp\!\left(Y^{\mathrm{NN}}\right).
\]
For a physical coarse interval \(K\), the corresponding predicted
conductance is then
\[
    c_K^{\mathrm{NN}}
    =
    \frac{\kappa_{\mathrm{eff}}^{\mathrm{NN}}}{H_K},
\]
and the predicted homogeneous DtN matrix is
\[
    S_K^{\mathrm{NN}}
    =
    c_K^{\mathrm{NN}}
    \begin{bmatrix}
        1 & -1\\
        -1 & 1
    \end{bmatrix}.
\]
Thus, the neural network learns only the coefficient-dependent scalar
quantity on the reference interval, while the geometric scaling and the
algebraic structure of the DtN matrix are imposed exactly.

For each coefficient resolution \(N_\kappa=16,32,64\) and each contrast
regime, we generate \(1200\) independent local coefficient realizations.
The first \(1000\) samples are used for training and the remaining \(200\)
for validation, and a separate network is trained for each value of
\(N_\kappa\). For the low-contrast case, we use a fully connected network
with three hidden layers of width \(64\), while for the high-contrast case
we use four hidden layers of width \(96\). Both networks use the
\texttt{Tanh} activation function, a batch size of \(64\), and the Adam
optimizer~\cite{KingmaBa2015} with learning rate \(10^{-3}\). The low- and
high-contrast networks are trained for \(500\) and \(700\) epochs,
respectively.

For the global verification, let $\Omega=(0,1)$ be partitioned into
$N_c$ uniform coarse intervals, with
$
    N_c=5,10,20,40.
$
Each coarse interval contains $N_\kappa$ coefficient cells.  For every
combination of $N_\kappa$, $N_c$, and contrast, the reported errors are
averaged over $25$ independent global coefficient realizations.

In one dimension, the boundary of each coarse interval consists only
of its two endpoints.  Hence the coarse trace contains the complete
HDG trace information on the coarse skeleton.  Standard MsHDG is
therefore simply a blockwise static condensation of the fine-scale HDG
system, and $u_H^{\mathrm{MS}}$ and $u_h^{\mathrm{REF}}$ agree up to
numerical roundoff.  In the present experiments,
$E_{\mathrm{MS}}$ ranges approximately from $10^{-14}$ to $10^{-9}$.
Accordingly,
\[
    E_{\mathrm{NN}}\approx E_{\mathrm{model}},
\]
and $E_{\mathrm{model}}$ provides a direct measure of the error
introduced by replacing the exact local DtN operators with their neural
approximations.

Table~\ref{tab:1d-model-error} reports the mean modeling errors.  The
errors remain below $10^{-2}$ for all coefficient resolutions, coarse
partitions, and both contrast regimes.

\begin{table}[ht]
\centering
\caption{One-dimensional modeling error $E_{\mathrm{model}}$.  The
reported values are sample means over $25$ independent global
coefficient realizations.}
\label{tab:1d-model-error}
\begin{tabular}{cccccc}
\hline
Contrast & $N_\kappa$
& $N_c=5$ & $N_c=10$ & $N_c=20$ & $N_c=40$\\
\hline
$\{1,10\}$
& 16
& $1.28\times10^{-3}$
& $1.22\times10^{-3}$
& $7.15\times10^{-4}$
& $7.21\times10^{-4}$\\
& 32
& $5.25\times10^{-3}$
& $3.89\times10^{-3}$
& $2.96\times10^{-3}$
& $2.45\times10^{-3}$\\
& 64
& $6.21\times10^{-3}$
& $3.59\times10^{-3}$
& $3.12\times10^{-3}$
& $2.01\times10^{-3}$\\
\hline
$\{1,10^4\}$
& 16
& $1.54\times10^{-3}$
& $1.93\times10^{-3}$
& $8.90\times10^{-4}$
& $8.01\times10^{-4}$\\
& 32
& $4.20\times10^{-3}$
& $3.99\times10^{-3}$
& $2.97\times10^{-3}$
& $2.04\times10^{-3}$\\
& 64
& $8.42\times10^{-3}$
& $6.38\times10^{-3}$
& $3.81\times10^{-3}$
& $2.81\times10^{-3}$\\
\hline
\end{tabular}
\end{table}

As expected, increasing $N_\kappa$ generally makes the local learning
problem more difficult when the training-data budget is held fixed.
Nevertheless, the resulting global modeling error remains below $1\%$
for every case considered.  Moreover, the error does not increase as
the number of coarse blocks grows; in most cases it decreases as
$N_c$ is increased.  Thus, in this analytically transparent example,
the local neural approximation errors do not accumulate destructively
through the global MsHDG assembly.

It is also noteworthy that the neural approximation remains accurate
for the high-contrast coefficient set $\{1,10^4\}$.  Hence coefficient
contrast alone does not explain the difficulties observed later in the
two-dimensional high-contrast experiments.  The one-dimensional DtN
map is determined by a single effective conductance and therefore does
not contain the richer interface behavior, coefficient-dependent trace
directions, and connectivity effects that arise in two dimensions.

No computational acceleration is expected in this one-dimensional
example.  The exact local DtN map has the closed-form scalar
representation~\eqref{eq:one-dimensional-conductance}, and the global reduced system is
extremely small.  Indeed, the measured total online speedup is
approximately $0.98$--$1.00$ over the tested cases.  The
one-dimensional experiment should therefore be viewed as a verification
of the learned local DtN operator and the global assembly procedure,
rather than as a demonstration of computational speedup.  The efficiency
benefit of NN-MsHDG is examined in the two-dimensional experiments,
where construction of each coarse-block DtN operator requires
substantially larger fine-scale HDG solves.

\subsection{Two-dimensional lower-contrast experiments}
\label{subsec:two-dimensional-low-contrast}

We next consider $\Omega=[0,1]^2$ with
$
    \kappa\in\{1,10\}.
$
The domain is partitioned into a uniform $5\times5$ array of square coarse blocks with side length $L = 1/5$.  Each coarse block is subdivided into $32\times32$ fine squares, and each square is divided into two triangles.  Hence the complete fine mesh has mesh size
\[
    h=\frac{L}{32}=\frac{1}{160}.
\]
We test the three trace levels $n=0,1,2$ introduced in the preceding section.  After imposing the homogeneous boundary condition, the corresponding numbers of globally coupled coarse trace unknowns are $80$, $160$, and $320$, respectively.

For the trace-level comparison, we use the same $20$ global coefficient seeds for $n=0,1,2$.  Reusing the realizations across trace levels gives a paired comparison in which changes can be attributed to trace-space enrichment rather than to different permeability draws. Table~\ref{tab:low-contrast-accuracy} reports the sample mean and sample standard deviation of the three relative errors.  The mean measures the typical accuracy, while the standard deviation records the sensitivity of that accuracy to the heterogeneous coefficient realization.

\begin{table}[h]
    \centering
    \caption{Lower-contrast two-dimensional accuracy.  Each entry is the sample mean $\pm$ sample standard deviation over the same $20$ coefficient realizations used at every trace level.}
    \label{tab:low-contrast-accuracy}
    \small
    \begin{tabular}{ccccc}
        \toprule
        $n$
        & Global trace DoFs
        & $E_{\mathrm{MS}}$
        & $E_{\mathrm{NN}}$
        & $E_{\mathrm{model}}$ \\
        \midrule
        0 & 80
          & $0.09782\pm0.00918$
          & $0.09637\pm0.01069$
          & $0.01351\pm0.00362$ \\
        1 & 160
          & $0.06683\pm0.00513$
          & $0.06534\pm0.00840$
          & $0.01338\pm0.00446$ \\
        2 & 320
          & $0.03525\pm0.00315$
          & $0.04516\pm0.01394$
          & $0.03001\pm0.00964$ \\
        \bottomrule
    \end{tabular}
\end{table}

The standard MsHDG error decreases from approximately $9.8\%$ to $3.5\%$ as the trace space is enriched.  For $n=0$ and $n=1$, the NN-MsHDG and standard MsHDG errors relative to the fine-scale reference are nearly identical, and the direct modeling error is approximately $1.3\%$.  At $n=2$, the neural modeling error increases to approximately $3.0\%$, but the accelerated solution remains at the same overall accuracy level as the standard MsHDG solution.  The increase is consistent with the greater complexity of the learned output: from $n=1$ to $n=2$, the local trace dimension doubles from $16$ to $32$, while the total number of predicted entries increases from $152$ to $560$.

\begin{figure}[htbp]
	\centering
	\includegraphics[width=\textwidth]{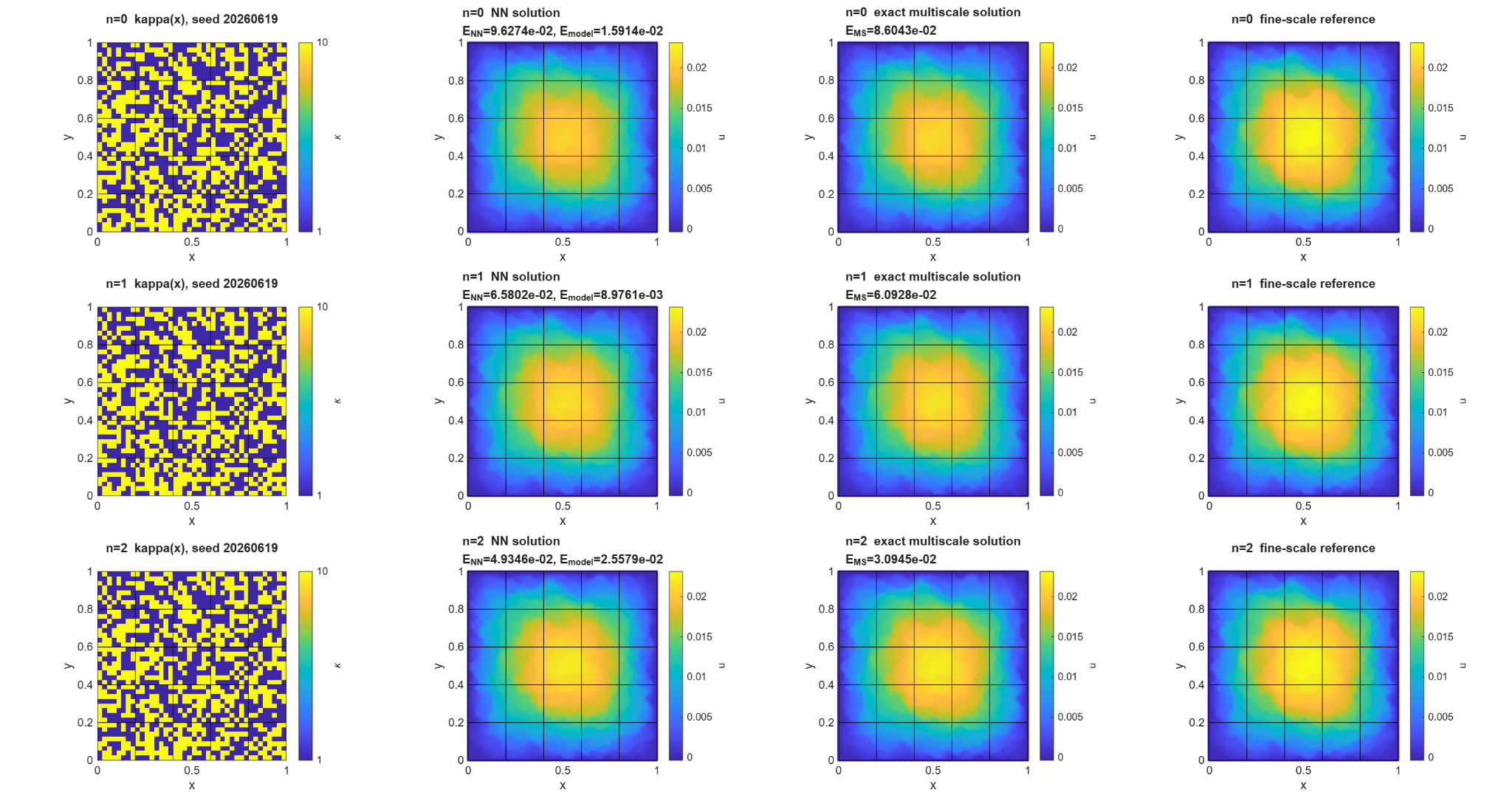}
	\caption{Representative lower-contrast comparison for the same
permeability realization at trace levels \(n=0,1,2\).
The columns show the permeability field, the NN-MsHDG solution, the
standard MsHDG solution, and the fine-scale HDG reference solution.}
	\label{fig:kappa10_same_seed_all_scales}
\end{figure}

\subsubsection{Online computational cost}
\label{subsec:low-contrast-timing}

The timing comparison is carried out in the lower-contrast regime because this is the regime in which the neural surrogate provides a sufficiently accurate replacement for the local HDG solvers.  The offline costs of data generation and network training are excluded.  Both MsHDG and NN-MsHDG timings were obtained on the same workstation described in Subsection~\ref{subsec:computational-environment}; the MATLAB parallel pool was disabled, coarse-block processing was serial, and MATLAB default internal multithreading was left unchanged.  Table~\ref{tab:low-contrast-timing} reports the mean wall-clock times, in seconds, from the final matched timing study.  That study used 15 repeats at each trace level and contrast, with the same seed used for the low- and high-contrast member of each matched pair.

\begin{table}[!htbp]
	\centering
	\caption{Online computational cost for the lower-contrast experiments.}
	\label{tab:low-contrast-timing}
	
	\small
	\setlength{\tabcolsep}{4pt}
	\renewcommand{\arraystretch}{1.15}
	
	\begin{tabular}{c c c c c c c c}
		\toprule
		$n$
		&
		\begin{tabular}[c]{@{}c@{}}
			Global trace\\
			DoFs
		\end{tabular}
		&
		$T_{\mathrm{asm}}^{\mathrm{NN}}$
		&
		$T_{\mathrm{asm}}^{\mathrm{MS}}$
		&
		$\mathcal S_{\mathrm{asm}}$
		&
		$T_{\mathrm{online}}^{\mathrm{NN}}$
		&
		$T_{\mathrm{online}}^{\mathrm{MS}}$
		&
		$\mathcal S_{\mathrm{online}}$
		\\
		\midrule
		0 & 80  & 6.333 & 53.180  & 8.397  & 12.270 & 59.190  & 4.824 \\
		1 & 160 & 6.208 & 98.150  & 15.811 & 11.979 & 103.852 & 8.669 \\
		2 & 320 & 6.320 & 192.947 & 30.530 & 12.160 & 198.703 & 16.341 \\
		\bottomrule
	\end{tabular}
\end{table}

The NN assembly time remains nearly unchanged as the trace space is enriched, whereas the standard MsHDG assembly time grows from approximately $53$ seconds to $193$ seconds because an increasing number of local fine-scale HDG problems must be solved.  Consequently, the assembly speedup increases approximately from $8.4$ to $30.5$.  The total online speedup increases from $4.8$ at $n=0$ to $16.3$ at $n=2$.  These results demonstrate the principal computational advantage of the proposed method: after training, inexpensive network evaluations replace repeated local fine-scale solves while the standard MsHDG global coupling is retained.

The timing results also confirm the computational bottleneck identified
in Section~2.  For the standard MsHDG method, the local-operator
construction and global assembly account for approximately \(90\%\),
\(95\%\), and \(97\%\) of the total online time for \(n=0,1,2\),
respectively.  Thus, most of the computational cost is not associated
with solving the global coarse skeleton system or reconstructing the
interior solution, but with the repeated fine-scale local solves required
to construct the coarse-block DtN operators.  This cost becomes more
pronounced as the trace space is enriched, since the local inverse must
be applied to an increasing number of coarse trace basis inputs.  The
advantage of replacing this stage by a neural surrogate is therefore
expected to be even greater for more complicated PDE systems, larger
local discretizations, or problems in which the local operators must be
reconstructed repeatedly.

\subsection{Two-dimensional high-contrast experiments}

We next consider the high-contrast regime
$
    \kappa\in\{1,10^4\}.
$
The coarse partition, fine mesh, coefficient representation, and trace spaces are
the same as in the lower-contrast experiments. Within this contrast regime, the
same $20$ held-out global coefficient realizations are used for the three trace
levels $n=0,1,2$.

Table~\ref{tab:high-contrast-accuracy} reports the global solution errors for the
fine-scale HDG reference solution $u_h^{\mathrm{REF}}$, the standard MsHDG
solution $u_H^{\mathrm{MS}}$, and the neural-network-based solution
$u_H^{\mathrm{NN}}$.

\begin{table}[ht]
\centering
\caption{High-contrast two-dimensional accuracy. Each entry is the sample
mean $\pm$ sample standard deviation over the same $20$ held-out coefficient
realizations used at every trace level.}
\label{tab:high-contrast-accuracy}
\begin{tabular}{ccccc}
\hline
$n$ & Global trace DoFs & $E_{\mathrm{MS}}$ & $E_{\mathrm{NN}}$
& $E_{\mathrm{model}}$\\
\hline
0 & 80
& $0.94864 \pm 0.01587$
& $0.96005 \pm 0.01058$
& $0.27945 \pm 0.22165$\\
1 & 160
& $0.84478 \pm 0.10386$
& $0.95584 \pm 0.01021$
& $0.69909 \pm 0.17307$\\
2 & 320
& $0.43394 \pm 0.18645$
& $0.94442 \pm 0.01340$
& $0.90271 \pm 0.04154$\\
\hline
\end{tabular}
\end{table}

An important observation is that the standard MsHDG approximation itself is
not sufficiently accurate in this high-contrast regime with the present
piecewise-linear coarse trace spaces. Although enrichment from $n=0$ to $n=2$
reduces $E_{\mathrm{MS}}$ substantially, the error remains approximately $43\%$
at the finest trace level considered. Thus, unlike the lower-contrast
experiment, the present choice of $M_H$ does not provide an accurate
coarse-scale approximation of the fine-grid solution for
$\kappa_{\mathrm{high}}=10^4$.

This behavior is consistent with the fact that high-contrast multiscale
problems generally require richer coarse spaces capable of representing
coefficient-dependent interface behavior. One possibility is to further
increase the dimension of $M_H$. A more systematic alternative is to use
coefficient-adapted spectral trace spaces, as in the spectral MsHDG
framework~\cite{EfendievLazarovMoonShi2014}. Investigating such choices is
beyond the scope of the present work, whose primary objective is to study
whether the fine-scale local MsHDG operators can be replaced efficiently by
neural-network surrogates.

For this reason, the errors relative to the fine-scale reference solution
should not be used alone to assess the quality of the neural local solver in
the high-contrast case. We instead compare the neural-network-predicted local
DtN operator directly with the corresponding operator produced by the
fine-scale HDG local solver. For a trace vector $\lambda$, we define the
relative local action error
\begin{equation}
    E_{\mathrm{act}}(\lambda;\alpha)
    =
    \frac{
    \left\|
       \left[S^{\mathrm{NN}}(\alpha)-S(\alpha)\right]\lambda
    \right\|_2}
    {\left\|S(\alpha)\lambda\right\|_2}.
    \label{eq:local-action-error}
\end{equation}
Here $S(\alpha)$ is the DtN matrix computed by the standard local HDG solver,
while $S^{\mathrm{NN}}(\alpha)$ is its neural-network approximation. Thus,
this quantity isolates the error introduced by the learned local solver and
does not include the approximation error associated with the choice of the
coarse trace space $M_H$.

Table~\ref{tab:dtn-action-errors} reports the mean relative action errors on
held-out reference-block coefficient realizations. The smooth-action
diagnostic uses smoothly varying trace probes, whereas the random-action
diagnostic uses randomly generated constant-free probes. These quantities are
discrete coefficient-vector errors and are used only as diagnostics of the
local DtN approximation.

\begin{table}[ht]
\centering
\caption{Relative local DtN action errors on held-out reference-block data.}
\label{tab:dtn-action-errors}
\begin{tabular}{ccccc}
\hline
& \multicolumn{2}{c}{$\kappa\in\{1,10\}$}
& \multicolumn{2}{c}{$\kappa\in\{1,10^4\}$}\\
\cline{2-5}
$n$
& Smooth action & Random action
& Smooth action & Random action\\
\hline
0 & $6.70\times10^{-4}$ & $6.88\times10^{-4}$
  & $2.00\times10^{-2}$ & $2.33\times10^{-2}$\\
1 & $6.39\times10^{-4}$ & $8.00\times10^{-4}$
  & $9.74\times10^{-3}$ & $1.06\times10^{-2}$\\
2 & $8.02\times10^{-3}$ & $1.02\times10^{-2}$
  & $7.96\times10^{-3}$ & $8.99\times10^{-3}$\\
\hline
\end{tabular}
\end{table}

For $n=0$ and $n=1$, the local DtN action errors are larger in the
high-contrast regime than in the lower-contrast case, indicating that the
coefficient-to-operator map becomes more difficult to approximate as the
contrast increases. Nevertheless, the errors remain on the order of
$10^{-2}$ on the tested probe families. At $n=2$, the local action errors for
the two contrast regimes are of comparable magnitude.

The contrast between these relatively small local operator-action errors and
the much larger global modeling error $E_{\mathrm{model}}$ in
Table~\ref{tab:high-contrast-accuracy} indicates that the high-contrast
problem is sensitive to perturbations of the local DtN operators. In
particular, average accuracy on the present probe families does not by itself
guarantee accuracy of the assembled global solution. This may reflect the
conditioning of the global skeleton system or the presence of
coefficient-dependent trace directions that are especially important in the
high-contrast regime. A systematic analysis of this error amplification,
together with the use of coefficient-adapted coarse trace spaces, is left for
future work.

For completeness, Table~\ref{tab:high-contrast-timing} reports the online
computational cost. The same timing protocol as in the lower-contrast
experiment is used.

\begin{table}[ht]
\centering
\caption{Online computational cost for the high-contrast experiments. The
reported values are mean wall-clock times in seconds over the $15$ repeats
of the matched timing study.}
\label{tab:high-contrast-timing}
\begin{tabular}{cccccccc}
\hline
$n$ & Global trace DoFs
& $T_{\mathrm{asm}}^{\mathrm{NN}}$
& $T_{\mathrm{asm}}^{\mathrm{MS}}$
& $S_{\mathrm{asm}}$
& $T_{\mathrm{online}}^{\mathrm{NN}}$
& $T_{\mathrm{online}}^{\mathrm{MS}}$
& $S_{\mathrm{online}}$\\
\hline
0 & 80
& 6.358 & 53.426 & 8.403
& 12.370 & 59.436 & 4.805\\
1 & 160
& 6.198 & 97.901 & 15.794
& 11.952 & 103.624 & 8.670\\
2 & 320
& 6.625 & 203.166 & 30.667
& 12.461 & 208.889 & 16.764\\
\hline
\end{tabular}
\end{table}

As in the lower-contrast experiment, replacing the local fine-scale solves by
network evaluations substantially reduces the online computational cost. The
assembly speedup increases from approximately $8.4$ to $30.7$, and the total
online speedup increases from approximately $4.8$ to $16.8$ as the trace
space is enriched. These timing results confirm that the computational
advantage of the neural local solver persists at high contrast. The principal
remaining questions are therefore approximation issues: the selection of a
coarse trace space suitable for high-contrast media and the sensitivity of
the assembled skeleton problem to errors in the learned local DtN operators.

\section{Concluding Remarks}\label{sec:conclusion}

In this work, we developed a neural-network-accelerated multiscale
hybridizable discontinuous Galerkin method for elliptic problems with
heterogeneous coefficients.  The proposed method preserves the local-to-global structure of the standard MsHDG formulation. On each coarse block, a fine-scale HDG solver defines a discrete Dirichlet-to-Neumann operator that maps the coarse trace data to the corresponding numerical flux response. These local operators are then assembled through the standard MsHDG global skeleton equations. The main idea of NN-MsHDG is to replace the repeated construction of the local HDG operators by a neural-network surrogate trained on coefficient fields defined on a reference block, while leaving the global coupling and local reconstruction procedures unchanged.

The numerical experiments demonstrate both the potential and the
current limitations of this approach.  The one-dimensional
fixed-sample study varies the local coefficient resolution
$N_{\kappa}=16,32,64$ while keeping 1000 training and 200 validation samples for each coefficient resolution and contrast regime.
  Increasing $N_{\kappa}$ makes the local
conductance-learning problem more demanding, but the assembled error
does not accumulate destructively: at $N_c=40$, the mean
$E_{\mathrm{model}}$ is below $3\times10^{-3}$ for every tested
coefficient resolution and for both contrast regimes.  In the
two-dimensional lower-contrast regime
\(\kappa\in\{1,10\}\), the NN-MsHDG solutions remain close to the
standard MsHDG solutions.  The modeling errors are approximately
\(1.3\%\) for trace levels \(n=0,1\) and approximately \(3.0\%\) for
\(n=2\).  At the same time, replacing the local fine-scale solves by
neural-network evaluations substantially reduces the online
computational cost.  The total online speedup increases from
approximately \(4.8\) to \(16.3\) as the coarse trace space is enriched.
These results show that the proposed method can provide meaningful
acceleration when the coarse-block local systems are sufficiently large
and the learned DtN operators are accurate in the directions relevant
to the assembled solution.

The high-contrast experiments
\(\kappa\in\{1,10^4\}\) identify an important limitation of the present
network and training strategy.  In this regime, the NN-MsHDG solution
does not remain close to the standard MsHDG solution, and the
deterioration becomes more pronounced as the coarse trace space is
enriched.  One possible explanation is that high coefficient contrast
produces local and global operators with strongly separated scales and
interface directions that are highly sensitive to perturbations.
Errors in the learned local DtN matrices may therefore be amplified
during assembly and inversion of the global skeleton system.  The
coefficient-to-operator map may also become more difficult to
approximate because small changes in the topology of highly conductive
regions can produce substantial changes in the effective flux response.

The choice of the coarse trace space may provide another important
direction for improvement.  The present study uses fixed discontinuous
piecewise-linear trace spaces.  For high-contrast multiscale problems,
however, polynomial trace spaces may fail to represent coefficient-dependent interface behavior efficiently.  Spectral multiscale HDG
methods construct reduced trace spaces from local snapshot functions
and coefficient-dependent spectral decompositions, thereby incorporating
important features of the heterogeneous medium into the coarse skeleton
space; see \cite{EfendievLazarovMoonShi2014}.  Combining neural
approximation of local DtN operators with such coefficient-adapted
spectral trace spaces is therefore a promising direction for future
work.  This extension is nontrivial because the basis of the trace
space would itself depend on the local coefficient field, and the
learned operators would have to be represented in coefficient-dependent
coordinates.

A rigorous study of the high-contrast regime remains necessary.
In particular, future analysis should clarify the relationship between
local DtN approximation errors, conditioning of the assembled skeleton
operator, and the resulting global solution error.  It would also be
useful to identify operator norms and solution-relevant trace directions
whose control guarantees stability of the learned global system.
Possible computational improvements include structure-preserving
parameterizations of the DtN matrix, coefficient- or energy-scaled loss
functions, adaptive sampling of difficult coefficient configurations,
and trace probes generated from representative global solutions.

Finally, the proposed framework is not restricted to the Darcy problem
considered here.  It is applicable whenever the standard HDG mechanism
is available: local element or subdomain unknowns are eliminated
independently, trace variables are globally coupled on the mesh
skeleton, and the local solver defines an input-output map from trace
and source data to numerical fluxes.  This includes other elliptic and
diffusion problems and, after temporal discretization or linearization,
time-dependent and nonlinear models.  The potential benefit may be even
greater in such applications because the local operators must be
reconstructed repeatedly at different time steps, parameter values, or
nonlinear iterations.


\begin{thebibliography}{10}
 


\bibitem{BoutilierBrennerMiguez2025}
M.~Boutilier, K.~Brenner, and L.~Miguez,
\emph{Learning local Dirichlet-to-Neumann maps of nonlinear elliptic
PDEs with rough coefficients},
ESAIM Proc. Surveys \textbf{81} (2025), 16--32.



\bibitem{ChenDingLiWright2024}
S.~Chen, Z.~Ding, Q.~Li, and S.~J. Wright,
\emph{A reduced order Schwarz method for nonlinear multiscale elliptic
equations based on two-layer neural networks},
J. Comput. Math. \textbf{42} (2024), no.~2, 570--596.

\bibitem{CockburnGopalakrishnanLazarov2009}
B.~Cockburn, J.~Gopalakrishnan, and R.~Lazarov,
\emph{Unified hybridization of discontinuous Galerkin, mixed, and
continuous Galerkin methods for second order elliptic problems},
SIAM J. Numer. Anal. \textbf{47} (2009), no.~2, 1319--1365.

\bibitem{CockburnGopalakrishnanSayas2010}
B.~Cockburn, J.~Gopalakrishnan, and F.-J.~Sayas,
\emph{A projection-based error analysis of HDG methods},
Math. Comp. \textbf{79} (2010), no.~271, 1351--1367.

\bibitem{CockburnQiuShi2012}
B.~Cockburn, W.~Qiu, and K.~Shi,
\emph{Conditions for superconvergence of HDG methods for second-order
elliptic problems},
Math. Comp. \textbf{81} (2012), no.~279, 1327--1353.




\bibitem{DuStechmann2024}
S.~Du and S.~N.~Stechmann,
\emph{Element learning: A systematic approach of accelerating finite
element-type methods via machine learning, with applications to radiative
transfer},
J. Comput. Math. \textbf{44} (2026), no.~1, 1--34.



\bibitem{EHanJentzen2017}
W.~E, J.~Han, and A.~Jentzen,
\emph{Deep learning-based numerical methods for high-dimensional
parabolic partial differential equations and backward stochastic
differential equations},
Commun. Math. Stat. \textbf{5} (2017), no.~4, 349--380.

\bibitem{EfendievLazarovMoonShi2014}
Y.~Efendiev, R.~Lazarov, M.~Moon, and K.~Shi,
\emph{A spectral multiscale hybridizable discontinuous Galerkin method
for second order elliptic problems},
Comput. Methods Appl. Mech. Engrg. \textbf{292} (2015), 243--256.

\bibitem{EfendievLazarovShi2015}
Y.~Efendiev, R.~Lazarov, and K.~Shi,
\emph{A multiscale HDG method for second order elliptic equations.
Part I: Polynomial and homogenization-based multiscale spaces},
SIAM J. Numer. Anal. \textbf{53} (2015), no.~1, 342--369.



\bibitem{HanJentzenE2018}
J.~Han, A.~Jentzen, and W.~E,
\emph{Solving high-dimensional partial differential equations using
deep learning},
Proc. Natl. Acad. Sci. USA \textbf{115} (2018), no.~34, 8505--8510.

\bibitem{HendrycksGimpel2016}
D.~Hendrycks and K.~Gimpel,
\emph{Gaussian error linear units ({GELU}s)},
arXiv preprint arXiv:1606.08415, 2016.



\bibitem{KabariaLewCockburn2015}
H.~Kabaria, A.~J. Lew, and B.~Cockburn,
\emph{A hybridizable discontinuous Galerkin formulation for non-linear
elasticity},
Comput. Methods Appl. Mech. Engrg. \textbf{283} (2015), 303--329.

\bibitem{KhooLuYing2021}
Y.~Khoo, J.~Lu, and L.~Ying,
\emph{Solving parametric PDE problems with artificial neural networks},
European J. Appl. Math. \textbf{32} (2021), no.~3, 421--435.

\bibitem{KingmaBa2015}
D.~P. Kingma and J.~Ba,
\emph{Adam: A method for stochastic optimization},
in: Proceedings of the 3rd International Conference on Learning
Representations (ICLR), 2015.

\bibitem{KochkovSmithAlievaWangBrennerHoyer2021}
D.~Kochkov, J.~A. Smith, A.~Alieva, Q.~Wang, M.~P. Brenner, and
S.~Hoyer,
\emph{Machine learning-accelerated computational fluid dynamics},
Proc. Natl. Acad. Sci. USA \textbf{118} (2021), no.~21, e2101784118.

\bibitem{KovachkiLiLiuAzizzadenesheliBhattacharyaStuartAnandkumar2023}
N.~Kovachki, Z.~Li, B.~Liu, K.~Azizzadenesheli, K.~Bhattacharya,
A.~Stuart, and A.~Anandkumar,
\emph{Neural operator: Learning maps between function spaces with
applications to PDEs},
J. Mach. Learn. Res. \textbf{24} (2023), no.~89, 1--97.




\bibitem{LiKovachkiAzizzadenesheliLiuBhattacharyaStuartAnandkumar2021}
Z.~Li, N.~Kovachki, K.~Azizzadenesheli, B.~Liu, K.~Bhattacharya,
A.~Stuart, and A.~Anandkumar,
\emph{Fourier neural operator for parametric partial differential
equations},
International Conference on Learning Representations, 2021.

\bibitem{LiKovachkiChoyEtAl2023}
Z.~Li, N.~Kovachki, C.~Choy, B.~Li, J.~Kossaifi, S.~Otta,
M.~A.~Nabian, M.~Stadler, C.~Hundt, K.~Azizzadenesheli,
and A.~Anandkumar,
\emph{Geometry-informed neural operator for large-scale 3D PDEs},
Adv. Neural Inf. Process. Syst. \textbf{36} (2023), 35836--35854.

\bibitem{LoshchilovHutter2017}
I.~Loshchilov and F.~Hutter,
\emph{SGDR: Stochastic gradient descent with warm restarts},
in: Proceedings of the 5th International Conference on Learning
Representations (ICLR), 2017.

\bibitem{LoshchilovHutter2019}
I.~Loshchilov and F.~Hutter,
\emph{Decoupled weight decay regularization},
in: Proceedings of the 7th International Conference on Learning
Representations (ICLR), 2019.

\bibitem{LuJinPangZhangKarniadakis2021}
L.~Lu, P.~Jin, G.~Pang, Z.~Zhang, and G.~E. Karniadakis,
\emph{Learning nonlinear operators via {DeepONet} based on the universal
approximation theorem of operators},
Nature Mach. Intell. \textbf{3} (2021), no.~3, 218--229.

\bibitem{LuWangXu2022}
Y.~Lu, L.~Wang, and W.~Xu,
\emph{Solving multiscale steady radiative transfer equation using
neural networks with uniform stability},
Res. Math. Sci. \textbf{9} (2022), no.~45.



\bibitem{MaFu2021}
W.~Ma and S.~Fu,
\emph{A hybridizable discontinuous Galerkin generalized multiscale
finite element method for highly heterogeneous linear elasticity
problems},
J. Comput. Appl. Math. \textbf{383} (2021), 113124.



\bibitem{NguyenPeraireCockburn2009Nonlinear}
N.~C. Nguyen, J.~Peraire, and B.~Cockburn,
\emph{An implicit high-order hybridizable discontinuous Galerkin method
for nonlinear convection--diffusion equations},
J. Comput. Phys. \textbf{228} (2009), no.~23, 8841--8855.

\bibitem{NguyenPeraireCockburn2010}
N.~C. Nguyen, J.~Peraire, and B.~Cockburn,
\emph{A hybridizable discontinuous Galerkin method for Stokes flow},
Comput. Methods Appl. Mech. Engrg. \textbf{199} (2010),
no.~9--12, 582--597.

\bibitem{NguyenPeraireCockburn2011NavierStokes}
N.~C. Nguyen, J.~Peraire, and B.~Cockburn,
\emph{An implicit high-order hybridizable discontinuous Galerkin method
for the incompressible Navier--Stokes equations},
J. Comput. Phys. \textbf{230} (2011), no.~4, 1147--1170.

\bibitem{NguyenPeraireCockburn2011Maxwell}
N.~C. Nguyen, J.~Peraire, and B.~Cockburn,
\emph{Hybridizable discontinuous Galerkin methods for the time-harmonic
Maxwell's equations},
J. Comput. Phys. \textbf{230} (2011), no.~19, 7151--7175.



\bibitem{RaissiPerdikarisKarniadakis2019}
M.~Raissi, P.~Perdikaris, and G.~E. Karniadakis,
\emph{Physics-informed neural networks: A deep learning framework for
solving forward and inverse problems involving nonlinear partial
differential equations},
J. Comput. Phys. \textbf{378} (2019), 686--707.



\bibitem{SerranoLeBoudecKassaiKoupaiWangYinVittautGallinari2023}
L.~Serrano, L.~Le Boudec, A.~Kassa\"{\i} Koupa\"{\i},
T.~X.~Wang, Y.~Yin, J.-N.~Vittaut, and P.~Gallinari,
\emph{Operator learning with neural fields: Tackling PDEs on general
geometries},
Adv. Neural Inf. Process. Syst. \textbf{36} (2023), 70581--70611.



\bibitem{TanoRagusa2021}
M.~E. Tano and J.~C. Ragusa,
\emph{Sweep-Net: An artificial neural network for radiation transport
solves},
J. Comput. Phys. \textbf{426} (2021), 109757.



\bibitem{YangShiFu2019}
Y.~Yang, K.~Shi, and S.~Fu,
\emph{Multiscale hybridizable discontinuous Galerkin method for flow
simulations in highly heterogeneous media},
J. Sci. Comput. \textbf{81} (2019), 1712--1731.

\bibitem{Yu2018}
B.~Yu,
\emph{The deep Ritz method: A deep learning-based numerical algorithm
for solving variational problems},
Commun. Math. Stat. \textbf{6} (2018), no.~1, 1--12.

\end{thebibliography}
\end{document}